\documentstyle[12pt]{article}
\input epsf
\newlength{\rig}

\newlength{\hei}
\newcommand{\fgr}[2]{
\setlength{\rig}{0.07\textwidth}
\setlength{\hei}{7.8cm}
\begin{figure}
\rule{\rig}{0in}
\epsfxsize=13cm
\epsffile{#1}
\caption{}\label{#2}
\end{figure}
}

\newcommand{\func}[1]{{\rm #1} \,}

\newcommand{\limfunc}[1]{{\rm{#1}}}
\newcommand{\text}[1]{{\mbox {#1}}}
\newtheorem{theorem}{Theorem}

\newtheorem{proposition}{Proposition}
\newtheorem{remark}{Remark}
\newtheorem{corollary}{Corollary}
\newtheorem{definition}{Definition}
\newtheorem{example}{Example}
\begin{document}

\setcounter{page}{1}

\author{{\bf Anton Savin\footnotemark\rule{5pt}{0pt} and
\addtocounter{footnote}{-1}
Boris Sternin\footnotemark}%
\addtocounter{footnote}{-1}\thanks{
Supported by the Russian Foundation for Basic Research under
grants No.~97-01-00703 and 97-02-16722a, Arbeitsgruppe Partielle
Differentialgleichungen und Komplexe Analysis, Institut f\"ur Mathematik,
Universit\"at Potsdam, and the International Soros Foundation.}
\\[3mm]
Moscow State University\\[3mm]
e-mail: antonsavin@mtu-net.ru\\[3mm]
e-mail: sternin@math.uni-potsdam.de}
\title{\bf Elliptic Operators in Even Subspaces}
\date{}
\maketitle

\vspace{15mm}

\begin{abstract}
An elliptic theory is constructed for operators acting in subspaces
defined via even pseudodifferential projections. Index formulas are
obtained for operators on compact manifolds without boundary and for
general boundary value problems. A connection with Gilkey's theory of
$\eta$-invariants is established.
\end{abstract}

\vspace{15mm}


\vspace{15mm}

{\bf 1991 AMS classification}: Primary 58G03, Secondary 58G10, 58G12,
58G25, 19K56

\vfill

\newpage

\tableofcontents

\section*{Introduction}
\addcontentsline{toc}{section}{Introduction}

We study elliptic operators in subspaces
defined by pseudodifferential projections, more precisely, {\em even\/}
pseudodifferential projections. Such projections appear already in
classical boundary value problems. Indeed, when we reduce an elliptic
boundary value problem to the boundary, the corresponding operator
acts in subspaces defined by pseudodifferential
projections. Moreover, projections prove to be very helpful in the
construction of a Fredholm theory for elliptic operators violating the
well-known Atiyah--Bott condition \cite{AtBo2}. It is shown in \cite{ScSS18}
that such a theory can be constructed in subspaces defined as the ranges of
pseudodifferential projections in Sobolev spaces. These results actually
go back to the classical Hardy spaces in which (and only in which) valid
Fredholm theory for the Cauchy--Riemann operator can be constructed. The
class of even projections is an important class of pseudodifferential
projections. The precise definition will be given later, but for now we
point out that in the subspaces defined by such projections we not only
prove the finiteness theorem, but also present the corresponding index
formula.

The most essential and fundamental distinction of elliptic theory in
subspaces defined by pseudodifferential projections from a similar theory
in Sobolev spaces is as follows. Although the ellipticity condition, just
as in the classical case, is expressed in terms of the principal symbols
of the main operator and the projections, the index of an elliptic
operator in subspaces is determined by neither the principal nor
even the complete symbols. This results in the necessity to give some
numerical characteristic of pseudodifferential projections or, which is
the same, of the subspaces they define. One can obtain such a numerical
characteristic, playing in a certain sense the role of the dimension of a
projection (in the finite-dimensional case it is equal to the rank of the
projection), at least in the class of even projections. This notion is
undoubtedly fundamental in our theory. In these terms we obtain an index
formula in the situation of compact manifolds without boundary as well as
for general boundary value problems. The relation between the notion of
``dimension'' and the $\eta$-invariant of Atiyah--Patodi--Singer \cite{APS1}
is established.

The paper is organized as follows.

In the first section, we consider subspaces defined as the ranges of
pseudodifferential projections on an odd-dimesional compact manifold $M$
without boundary.\footnote{Let us note that all our constructions are
determined by the subspaces themselves and are independent of the choice
of projections onto them.} More precisely, we assume that the projection
is {\em even\/}, in the sense that its principal symbol is an even function
with respect to the cotangent variables.

It turns out that in the class of such subspaces there exists a uniquely
defined (up to a normalization) analog of the notion of {\em dimension\/}
of a finite-dimensional vector space (i.e., a homotopy invariant additive
functional). Every normalization is a choice of dimensions for spaces
of sections of vector bundles over $M$.

Even projections have the following property: the group of stable homotopy
classes of even projections is rationally generated by elements that
differ from the projections on spaces of sections of vector bundles on $M$
by finite rank operators. This statement is actually a consequence of the
fact that the obstruction to the stable homotopy of projections (modulo
finite rank operators) lies in the group $K\left( P^{*}M\right) /K\left(\!
M\!\right) $ (here $P^{*}\!M\!$ is the projectivization of the cosphere
bundle), which is a torsion group.

In the second section we study elliptic operators acting in subspaces
defined by pseudodifferential projections:
\[
D:H_1\longrightarrow H_2,\quad H_{1,2}=\limfunc{Im}\,P_{1,2}.
\]
Namely, we show that the index of elliptic operators of this form is
represented as a sum of two homotopy invariant terms, one of which is
determined by the principal symbol of the operator $D$ (and is a homotopy
invariant of the principal symbol of the problem), while the second term
is determined only by the subspaces where the operator $D$ acts and is a
homotopy invariant of the subspaces.

Let us note that there is no decomposition of this kind in the class of
{\em all } elliptic operators. It can be shown that here, as well as in the
theory of spectral boundary value problems (see \cite{SaSSc1,Sav1}), there
is an obstruction to such a decomposition, which is closely related to the
spectral flows of periodic families of operators (see Sec.~\ref{perv}).
That is why, to obtain an index formula, one must necessarily take
narrower classes of operators.

In Sections \ref{genbvp}, \ref{genspec}, and \ref{indevn} we give
an application of the theory to boundary value problems.
The general boundary value problems \cite{ScSS18} are
considered in Sec.~\ref{genbvp}. They have he form
\begin{equation}
\left\{
\begin{array}{l}
Du=f, \\
Bu=g\in {\mathrm{}Im}\,P,
\end{array}
\right.   \label{bvp}
\end{equation}
where $D$ is an elliptic differential operator on a smooth manifold $M$
with boundary $\partial M$, $B$ is a boundary operator, and $P$ is some
pseudodifferential projection on $\partial M$. This class of boundary
value problems, on the one hand, contains all classical boundary value
problems with the Lopatinskii condition (e.g., see \cite{Hor3}). On the
other hand, for any elliptic operator $D$ there exists a Fredholm boundary
value problem in this class. In this section all the necessary definitions
are given. Examples are presented.

Sections \ref{genspec} and \ref{indevn} deal with the index
computation for general boundary value problems. In
Sec.~\ref{genspec} the problem is reduced to a certain spectral boundary
value problem \cite{NScSS3}. In a collar neighborhood of the boundary
with the normal coordinate $t$, it has the form
\begin{equation}
\left\{
\begin{array}{l}
\left( \frac \partial {\partial t}+A\left( t\right) \right) u=f, \\
P\left. u\right| _{\partial M}=g\in {\mathrm{}Im}\,P,
\end{array}
\right.   \label{bvpspec}
\end{equation}
where $A\left( 0\right) $ is a pseudodifferential operator on the boundary,
whose homogeneous principal symbol on the bundle of cotangent spheres to the
boundary is equal to
\[
\sigma \left( A\left( 0\right) \right) =2\sigma \left( P\right) -1.
\]
The reduction is understood in the sense that the corresponding
Fredholm operators have the same indices.

In Sec.~\ref{indevn} the class of even boundary value problems on an
even-dimensional manifold $M$ is considered. These are boundary value problems
of the form (\ref{bvp}) with an even projection $P$. Under this condition,
we obtain the following index formula for the spectral boundary value
problem (\ref{bvpspec}):
\begin{equation}
{\mathrm{}ind}\left( D,P\right) =\frac 12\,{\mathrm{}ind}_t\left(
\sigma \left(\tilde D\right)
\right) -d\left( P\right) ,  \label{indmain}
\end{equation}
where $\sigma(\tilde D)$ is the elliptic symbol on the double of $M$
obtained by (continuously) gluing the symbols $\sigma(D) \left( \xi
\right) $ and $\sigma(D) \left( -\xi \right), $ and $d(P)$ is the
above-mentioned dimension functional for the trivial normalization. Let us
note that, by virtue of the reduction carried out in the first part of the
paper, formula (\ref{indmain}) solves the index problem for general
boundary value problems (\ref{bvp}) in the case of even projections $P.$

The third application of the introduced notion of $d$-dimension is related
to the $\eta$-invariant of Atiyah--Patodi--Singer \cite{APS1}. Namely, the
invariant $d$ of an even pseudodifferential projection is closely
connected with the theory of $\eta$-invariants of self-adjoint elliptic
operators of even order on odd-dimensional manifolds. More precisely, the
$\eta$-invariant of an {\em admissible\/} operator \cite{Gil7} (see also
Sec.~\ref{defadm}) in this case is equal to the introduced ``dimension''
of its nonnegative spectral subspace for the trivial normalization. It
follows that the fractional parts of the "dimension" and of the
$\eta$-invariant define a homomorphism
$$
K\left( P^{*}M\right) /K\left( M\right) \longrightarrow
{\bf Z}\left[\frac 12\right]
\,\func{mod}\,{\bf Z}.
$$
By virtue of this identification, the index formula (\ref{indmain})
can be interpreted as an analog of the Atiyah--Patodi--Singer formula
\cite{APS1}. As a corollary to the index formula in subspaces, we
obtain (see Corollary \ref{cor1}) a topological expression
for the fractional part of the $\eta$-invariant on the subgroup
$$
{\rm ker}\left\{\pi^*:
K\left( P^{*}M\right)/K(M) \longrightarrow K\left( S^{*}M\right)/K(M)
\right\},
\quad\pi:S^*M\to P^*M.
$$
While the index formula for even boundary value problems shows the
cobordism invariance of the fractional part of the doubled $\eta$-invariant
(see Corollaries \ref{cor2} and \ref{cor3}), it turns
out that in both cases the $\eta$-invariant has at most $2$ in
the denominator
or is even an integer. Nevertheless, the problem of the nontriviality of
this fractional part remains open.

In the last section of the paper we consider several examples.

The authors are grateful to Prof.~A.~S.~Mishchenko for a number of
valuable remarks he made when this work was reported at his seminar in
Moscow State University in fall 1998. He also suggested to include the
applications concerning the boundary value problems, which undoubtedly
improved the paper. We would like to thank V.~E.~Nazaikinskii for
constructive suggestions that helped us improve the original
version of the
paper. Finally, we are grateful to Prof.~P.~Gilkey, the discussion with
whom on the subject of this paper was extremely useful for us.

\section{Even pseudodifferential projections\label{perv}}

Let us consider the set $\limfunc{Proj}\,(M)$
of all pseudodifferential projections of order zero
acting in the spaces of smooth sections of vector bundles
on a closed manifold $M$. In this set it is impossible to
compare the ranks (the dimensions of the range) of projections
as one can do in the case of finite-dimensional spaces.
More precisely, on the space of pseudodifferential projections
there does not exist a homotopy invariant (with respect to
the operator $L^2$-norm) functional
\[
d:\limfunc{Proj}\,(M)\longrightarrow {\bf Z}\quad \left( \mbox{or even }%
{\bf R}\mbox{ or }{\bf C}\right) {\bf }
\]
that satisfies the (weak) additivity property
\[
d\left( P\oplus P^{\prime }\right) =d\left( P\right) +\dim \func{Im}%
P^{\prime },
\]
for an arbitrary projection $P$, where $P^{\prime }$ is a
finite-dimensional projection. Indeed, assuming the opposite, consider an
arbitrary periodic family $A_t$ of self-adjoint elliptic operators with a
nonzero spectral flow. Such families exist (e.g., see \cite{SaSSc1}).

Denote the corresponding family of projections on the nonnegative spectral
subspaces of the operators $A_t$ by $P_t$. From the definition of the
spectral flow $\limfunc{sf}A_t$ as the net number of eigenvalues of the
operators $A_t$ that pass through zero as the parameter $t$ varies, we
obtain
\begin{equation}
\limfunc{sf}A_t=d\left( P_1\right) -d\left( P_0\right)\neq 0. \label{dned}
\end{equation}
The family $A_t$ is, however, periodic. Consequently, $P_0=P_1$ and
$d\left( P_1\right) =d\left( P_0\right)$, which contradicts~(\ref{dned}).

In the remaining part of the section we consider the class of even
projections, where a similar functional can nevertheless be defined under
the condition that the manifold $M$ where the projections act is
odd-dimensional.

\begin{definition}
\label{def1}{\em A pseudodifferential projection }
\[
P:C^\infty \left( M,E\right) \longrightarrow C^\infty \left( M,E\right)
\]
{\em is called \ }even {\em (cf. \cite{Gil7}) if its
principal homogeneous symbol $\sigma(P)$ satisfies}
\begin{equation}
\sigma \left( P\right) \left( x,\xi \right) =\sigma \left( P\right) \left(
x,-\xi \right) \quad \mbox{\em for all  }\left( x,\xi \right) \in S^{*}M.
\label{A}
\end{equation}
\end{definition}
The set of all even pseudodifferential projections of order zero
is denoted by $\limfunc{Even}\left( M\right) .$ Let
\[
P^{*}M=\left. S^{*}M\right/ {\bf Z}_2
\]
be the bundle of projective spaces obtained as the quotient of the
cotangent sphere bundle $S^{*}M$ under the action of the antipodal
involution $\xi \rightarrow -\xi .$ The corresponding projection is
denoted by \ $\pi :S^{*}M\rightarrow P^{*}M.$ Then the symmetry condition
(\ref{A}) means that the homogeneous principal symbol $\sigma\left( P
\right) $ is the pullback of an endomorphism $\sigma^{\prime }\left(
P\right) $ over the projective bundle $P^{*}M :$
\begin{equation}
\sigma \left( P\right) =\pi ^{*}\sigma ^{\prime
}\left( P\right) ,\quad \sigma ^{\prime }\left( P\right):
\pi^*_PE\longrightarrow\pi^*_PE,
\label{bBB}
\end{equation}
where $\pi_P:P^*M\to M$ is the natural projection. Note that the
Eq.~(\ref{bBB}) implies that an even projection determines a vector bundle
on the projectivization $P^*M$,
\begin{equation}
\mbox{Im}\,\sigma'(P)\in
\limfunc{Vect}\left( P^{*}M\right).
\end{equation}

{Let }$P_{1,2}$ {be pseudodifferential projections
}
\[
P_{1,2}:C^\infty \left( M,E_{1,2}\right) \longrightarrow C^\infty \left(
M,E_{1,2}\right).
\]
Their {\em direct sum\/}$P_1\oplus
P_2${ \ is the projection}\
\[
P_1\oplus P_2=\left(
\begin{array}{ll}
P_1 & 0 \\
0 & P_2
\end{array}
\right) :
\begin{array}{c}
C^\infty \left( M,E_1\oplus E_2\right)
\end{array}
\longrightarrow
\begin{array}{c}
C^\infty \left( M,E_1\oplus E_2\right)
\end{array}
.
\]

Let us consider the following stable homotopy equivalence relation
on the set $\limfunc{Even}\left( M\right)$ of even projections.

\begin{definition}
{\em We say that two projections }$P_{1,2}\in \limfunc{%
Even}\left( M\right) $ {\em \ are \ }equivalent {\em and write}
{\em \ } $P_1\approx P_2${\em \ if for some even projection
}$P:C^\infty \left( M,F\right)
\longrightarrow C^\infty \left( M,F\right) $ {\em  there exists a homotopy
of even projections}
\[
P_1\oplus 0\oplus P\sim 0\oplus P_2\oplus P
\]
{\em as projections in the ambient space $C^\infty \left(
M,E_1\oplus E_2\oplus F\right) .$
}
\end{definition}

Now consider the Grothendieck group generated by the semigroup $\left.
\limfunc{Even}\left( M\right) \right/ \approx $ consisting of classes of
equivalent projections:
\begin{equation}
\label{omeg}
K\left( P_{ev}\left( M\right) \right)
\stackrel{{\rm def}}=K\left( \left. \limfunc{Even}\left(
M\right) \right/ \approx \right) .
\end{equation}
Each even projection $P$ defines an element in the group $K\left(
P_{ev}\left( M\right) \right) \ $, which we denote by $\left[ P\right] .$
It is clear that the mapping $q:K\left( P_{ev}\left( M\right) \right)
\longrightarrow K\left(P^{*}M\right)$ taking each even pseudodifferential
projection to the range of its principal symbol is a homomorphism of
abelian groups.

The above-mentioned properties of the rank of projections
are formalized in the following definition.

\begin{definition}
\label{dimen}
{\em
A group homomorphism}
\[
d:K\left( P_{ev}\left( M\right) \right) \longrightarrow {\bf R}
\]
{\em is called a}
dimension {\em %
if for any finite-dimensional projection }$P$ {\em \ we have}
\[
d\left( \left[ P\right] \right) ={\rm rank}\,P\equiv\dim \func{Im}P.
\]
\end{definition}

The following theorem describes all possible dimension functionals.

\begin{theorem}
\label{thmmn}
The dimensions of even subspaces
\[
d:K\left( P_{ev}\left( M\right) \right) \longrightarrow {\bf R}
\]
are in a one-to-one correspondence with homomorphisms
\[
\chi :K\left( M\right) \longrightarrow {\bf R},
\]
which will be called ``normalization homomorphisms.''
Moreover, for integer valued normalizations $%
\chi $ the dimension $d$ takes rational values whose denominators
can contain only powers of $2$:
\[
d:K\left( P_{ev}\left( M\right) \right) \longrightarrow {\bf Z}\left[
\frac 12\right] .
\]
\end{theorem}

\noindent {\em Proof\/}. Consider the sequence
\begin{equation}
0\longrightarrow {\bf Z}\stackrel{i}{\longrightarrow }
K\left( P_{ev}\left(
M\right) \right) \stackrel{q}{\longrightarrow }
K\left( P^{*}M\right) \longrightarrow
0,  \label{seq}
\end{equation}
where the map $i$ is induced by the homomorphism of
semigroups
$$
{\bf Z}_{+}\rightarrow
\limfunc{Even}\left( M\right) /\approx
$$
taking each nonnegative number
$k$ to the class of a projection of rank $k$ acting,
say, in the space $C^\infty \left(
M\right) .$ Let us prove that the sequence
(\ref{seq}) is exact.

First, we verify that $q$
is an epimorphism. Indeed, an arbitrary vector bundle
$\gamma \in \limfunc{Vect}\left( P^{*}M\right) $ can be
realized as a subbundle in some trivial bundle,
\[
\gamma \subset {\bf C}^N\in \limfunc{Vect}%
\left( P^{*}M\right).
\]
By lifting this embedding to the cotangent spheres, we obtain
\[
\pi ^{*}\gamma \subset {{\bf C}}^N\in \limfunc{Vect}%
\left( S^{*}M\right).
\]
It is obvious that the orthogonal projection
$\sigma =\sigma \left( x,\xi \right) $ on the subbundle
$\pi ^{*}\gamma $ is an even projection.
Consider an arbitrary pseudodifferential projection $P$ with
principal symbol $\sigma $\footnote{Following \cite{BiSo1},
we can define $P$ by the formula
$$
P=-\frac{1}{2\pi i}\int\limits_{|\lambda-1|=\varepsilon}
\left(\Pi-\lambda{\rm I}\right)^{-1}
d\lambda
$$
where $\Pi$ is an arbitrary pseudodifferential operator of order zero with
principal symbol $\sigma$, and the number $\varepsilon, 0<\varepsilon<1$,
is chosen in a way such that the circle $|\lambda-1|=\varepsilon$ contains
no eigenvalues of $\Pi$.}.
For the projection $P$ we obtain, by construction,
\[
q\left( \left[ P\right] \right) =\left[ \gamma \right] \in K\left(
P^{*}M\right) .
\]
This proves the exactness of the sequence in the third term.

The triviality of the composition $q\circ i$ is obvious. Let us check the
inclusion $\ker q\subset \func{Im}i.$ Suppose that
\begin{equation}
q\left( \left[ P\right] -\left[ 1_N\right] \right) =0,  \label{cond1}
\end{equation}
where $P$ is an even projection in $C^\infty \left( M,E\right),$ and $1_N$
is the identity in $C^\infty \left( M,{\bf C} ^N\right) $
(this does not restrict generality, since an arbitrary element of
$K\left( P_{ev}\left( M\right) \right) $ is representable in this form).
Condition (\ref{cond1}) means that (possibly, after adding a trivial
pair $\left( 1_{N^{\prime }},1_{N^{\prime }}\right) $ to the pair $\left(
P,1_N\right) $) we obtain an isomorphism
\[
\func{Im}\sigma' \left( P\right) \simeq {\bf C}^N
\]
of vector bundles over $P^*M$. It follows that these bundles are homotopic
as subbundles in the direct sum $\pi _P^{*}E\oplus {{\bf C}}^N$. Consider
an arbitrary homotopy joining them and denote the corresponding family of
projections by
\begin{eqnarray*}
\left\{ \sigma _t^{\prime }\right\} _{t=0,1} &:&\sigma _0^{\prime }=\sigma
^{\prime }\left( P\right) \oplus 0,\quad \sigma _1^{\prime }=0\oplus 1_N, \\
\func{Im}\sigma _t^{\prime } &\subset &\pi _P^{*}E\oplus {\bf C}^N\in
\limfunc{Vect}\left( P^{*}M\right) \quad \mbox{for all
}t\in \left[ 0,1\right] .
\end{eqnarray*}
The pullback of this homotopy to the cosphere bundle $S^{*}M$ will be
denoted by $\sigma _t.$ It follows from Statement 1 in \cite{MePi1} that
there exists a (continuous) covering homotopy of pseudodifferential
projections $P_t$,
\[
\sigma \left( P_t\right) =\sigma _t,
\]
such that $P_0=P\oplus 0$ and
$P_1$ differs from $0\oplus 1_N$ by a compact operator.
Accordingly, in the group $K(P_{ev}(M))$ we obtain
\[
\left[ P\right] -\left[ 1_N\right] =\left[ P_1\right] -\left[ 0\oplus
1_N\right].
\]
It can be shown (e.g., see \cite{Wojc1}) that projections
differing by a compact operator are homotopic
up to a finite rank projection. For example,
in the case of a positive relative index of projections\footnote{For
projections $P,Q$ with compact  difference,
the {\em relative index\/} is defined as the index of the
Fredholm operator $Q:\,{\rm Im}\,P\rightarrow {\rm Im}\,Q$:
\[
\limfunc{ind}\left( P,Q\right) \stackrel{{\rm def}}{=}
\limfunc{ind}\left( Q:\func{Im}P\rightarrow \func{Im}Q\right)
=-\limfunc{ind}\left( P:\func{Im}%
Q\rightarrow \func{Im}P\right).
\] }
$\limfunc{ind}\left( P_1,0\oplus 1_N\right) =n\geq 0$ we
obtain a homotopy of projections
\[
P_1\sim n\oplus 1_N,
\]
where $n$ is a rank $n$ projection in
$C^\infty \left( M,E\right) .$
For negative $\limfunc{ind}\left(
P_1,0\oplus 1_N\right) =n<0$, we have
\[
P_1\sim 0\oplus \left( 1_N-\left( -n\right) \right) ,
\]
where $-n$, as before, is a projection of
rank $-n>0.$ In the first case we find that
\[
\left[ P_1\right] -\left[ 0\oplus 1_N\right] =\left[ n\right] =i\left(
n\right)
\]
in the Grothendieck group $K\left( P_{ev}\left( M\right)\right) $,
and in the second case we also have
\[
\left[ P_1\right] -\left[ 0\oplus 1_N\right] =-\left[ -n\right] =i\left(
n\right) .
\]

Let us finally verify the exactness of the sequence
(\ref{seq}) in the first term. Suppose that for some $n>0$
we have
\[
i\left( n\right) =0.
\]
From the definition of the Grothendieck group it follows that
for some even projection $P\in \limfunc{Even}\left( M\right) $
there exists a homotopy of even projections
\[
n\oplus P\sim 0\oplus P.
\]
Moreover, without loss of generality, it can be assumed that the projection
$P$ is the unit operator: $P=1_N$. Let us denote this continuous homotopy
of even projections from $n\oplus P$ to $0\oplus P$ by $P_t.$ We obtain
\begin{equation}
\limfunc{ind}\left( P_0,P_1\right) =n\neq 0.  \label{rel1}
\end{equation}
Let us show that the fact that the projections
$P_t$ are even implies $\limfunc{ind}\left( P_0,P_1\right) =0.$
Indeed, without loss of generality it can be
assumed that the family
$\left\{ P_t\right\} $ consists of orthogonal projections
(the space of all projections can be linearly retracted to the
space of orthogonal projections,\footnote{The pseudodifferentiality
of the orthogonal projection
on the subspace ${\rm Im}P$ for an arbitrary pseudodifferential
projection $P$ follows from the identity
${\rm Im}P={\rm Im}PP^*$ and the self-adjointness of the operator
$PP^*$.}
and the relative index (\ref{rel1}) does not change under this
retraction).
Consider further an arbitrary
periodic family $\left\{ A_t\right\}_{t=0,1} $
of first-order self-adjoint pseudodifferential operators with
homogeneous principal symbols on the cotangent spheres equal to
$2\sigma \left( P_t\right) -1$
(so that the positive spectral projection of the
principal symbol of the operator
$A_t$ coincides with
$\sigma \left( P_t\right) $). The {\em spectral flow\/} of any
periodic family $\left\{ A_t\right\} _{t=0,1}$ can be expressed by the
cohomological formula
(\cite{SaSSc1}, cf. \cite{APS3}) %
\[
\limfunc{sf}\left\{ A_t\right\} _{t\in S^1}=\left\langle \limfunc{ch}\left(
\sigma _{+}\left( A_t\right) \right) \pi ^{*}\limfunc{Td}\left(
T^{*}M\otimes {\bf C}\right) ,\left[ S^1\times S^{*}M\right]
\right\rangle ,\quad \pi :S^1\times S^{*}M\rightarrow M.
\]
Here the vector bundle $\sigma _{+}\left( A_t\right) \in \limfunc{Vect}
\left( S^1\times S^{*}M\right) $ is generated by the nonnegative spectral
subspaces of the principal symbols $\sigma \left( A_t\right) $ (in our
case $\sigma _{+}\left( A_t\right) =\func{Im}\sigma \left( P_t\right) $),
and ${\rm Td}$ is Todd class of a vector bundle.

On the oriented manifold $S^{*}M$ the involution $\left( x,\xi \right)
\rightarrow \left( x,-\xi \right) $ reverses the orientation, while
the cohomology class
\[
\limfunc{ch}\left( \sigma _{+}\left( A_t\right) \right) \pi ^{*}\limfunc{Td}%
\left( T^{*}M\otimes {\bf C}\right)
\]
is invariant with respect to this involution. Thus, the spectral flow of
the family $\left\{ A_t\right\} _{t=0,1}$ is zero:
\[
\limfunc{sf}\left\{ A_t\right\} _{t=0,1}=0.
\]
Let us recall that the family of projections $P_t$ is a {\em generalized
spectral section\/} \cite{DaZh2} for the family $A_t.$ Then, by virtue of
one of the definitions of the spectral flow (see \cite{MePi1} or
\cite{DaZh2}), we obtain
\[
\limfunc{sf}\left\{ A_t\right\} _{t=0,1}=\limfunc{ind}\left( P_0,P_1\right)
.
\]
Hence,
\[
\limfunc{ind}\left( P_0,P_1\right) =0,
\]
which contradicts (\ref{rel1}). The exactness of the sequence
(\ref{seq}) is established.

Let us note that in terms of the sequence (\ref{seq}) the problem of
describing dimension homomorphisms $d$ is reduced to the problem of
closure of the following diagram to a commutative one:
\[
\begin{array}{ccc}
{{\bf Z}} & \stackrel{i}{\rightarrow } & K\left(
P_{ev}\left( M\right) \right) \\
\downarrow & \swarrow d &  \\
{\bf R} &  &
\end{array}
\]
(here ${\bf Z}\subset{\bf R}$ is the natural inclusion,
and the map $d$ is so far unknown). First, note that the groups
$K\left(
P_{ev}\left( M\right) \right) $ and $K\left( P^{*}M\right) $
contain subgroups generated by bundles on the base $M$:
\[
\begin{array}{ccccccc}
K\left( M\right) &\stackrel{}\rightarrow & K\left( P^{*}M\right),
& \qquad &
K\left(M\right) & \stackrel{\alpha}{\rightarrow}
       & K\left( P_{ev}\left( M\right) \right),\\
\left[ E\right] & \rightarrow & \left[ \pi _P^{*}E\right], &  & \left[
E\right] & \rightarrow & \left[ 1_{C^\infty \left( M,E\right) }\right].
\end{array}
\]
Both maps are monomorphisms
(this follows from the existence of a nonsingular vector field
on an odd-dimensional manifold $M$),
and moreover, these embeddings commute with
$q$. Thus, there is a partial splitting of the sequence (\ref{seq}):
\begin{equation}
\begin{array}{ccccccccc}
0 & \rightarrow & {{\bf Z}}\oplus K\left( M\right) & \stackrel{i\oplus\alpha
}{%
\rightarrow } & K\left( P_{ev}\left( M\right) \right) & \stackrel{q}{%
\rightarrow } & K\left( P^{*}M\right) /K\left( M\right) & \rightarrow & 0
\end{array}
.  \label{seq1}
\end{equation}
The quotient $K\left(P^{*}M\right) /K\left( M\right) $
is by \cite{Gil7} a purely
torsion group (and the torsion is only in powers of $2$).
By virtue of the sequence (\ref{seq1}), this implies that an arbitrary
dimension
\[
d:K\left( P_{ev}\left( M\right) \right) \longrightarrow {\bf R}
\]
is uniquely determined by its restriction to the subgroup
$K\left( M\right)$ (recall that on the first
term in the sum
${{\bf Z}}\oplus K\left( M\right)$ the dimension $d$ is already defined,
see Definition \ref{dimen}):
\[
\chi :K\left( M\right) \longrightarrow {\bf R.} %
\]
This completes the proof of the theorem.

\begin{remark}
\label{remmu}
{\em The dimension $d$ can be viewed as a generalization of the
relative index of projections, since}
\[
\limfunc{ind}\left( P_1,P_2\right) =d\left( \left[ P_1\right] -\left[
P_2\right] \right)
\]
{\em for two even projections $P_{1,2}$ with the same principal symbol.}
\end{remark}

For a more explicit expression for the dimension
$d\left( \left[ P\right]\right) $ with the trivial
normalization $\chi \equiv 0$ see Sec.~\ref{defadm}, where the
relation with Gilkey's $\eta$-invariants is established.

\section{Operators in subspaces. An index formula for even subspaces
\label{subsp1}}

\noindent
Consider two pseudodifferential projections
\[
P_{1,2}:C^\infty \left( M,E_{1,2}\right) \longrightarrow C^\infty \left(
M,E_{1,2}\right)
\]
of order zero on a manifold $M$
and an $m$th order pseudodifferential operator
\[
D:C^\infty \left( M,E_1\right) \longrightarrow C^\infty \left( M,E_2\right)
.
\]
Suppose that $D$ acts in the subspaces determined by the
projections, that is,
$$
D({\rm Im}\,P_1)\subset{\rm Im}\,P_2, \quad P_2DP_1=DP_1.
$$
Then the restriction
\[
D:\func{Im}P_1\longrightarrow \func{Im}P_2
\]
is called an {\em operator acting in subspaces\/}. There is a criterion
for the Fredholm property of operators of this form (see \cite{BiSo1} and
\cite{ScSS18}). Before stating it, let us introduce the notion of the
principal symbol of an operator in subspaces.
\begin{definition}
{\em The} principal symbol {\em of the operator}
\[
D:\func{Im}P_1\longrightarrow \func{Im}P_2,
\]
{\em acting in subspaces,
is the homomorphism of vector bundles over
$S^*M$ given by the restriction of the principal symbol of the operator}
$D$:
\[
\sigma \left( D\right) :\func{Im}\sigma \left( P_1\right) \longrightarrow
\func{Im}\sigma \left( P_2\right). \label{lbl1}
\]
\end{definition}

A symbol is called {\em elliptic\/} if it is an isomorphism. In this case
the operator $D$ is also called {\em elliptic\/}. In these terms, the
following statement is valid; the proof is essentially contained in
\cite{ScSS18}.
\begin{proposition}
An operator\
\[
D:H^s\left( E_1\right) \supset \func{Im}P_1\longrightarrow \func{Im}P_2
\subset H^{s-m}\left( E_2\right)
\]
is Fredholm if and only if it is elliptic.
\end{proposition}
The orders of the Sobolev spaces will be omitted for brevity in what follows.

Let us consider triples  $(D,P_1,P_2)$ representing
elliptic operators
\[
D:\func{Im}P_1\longrightarrow \func{Im}P_2
\]
such that $P_{1,2}$ are even projections in the sense of Definition 1.
These triples form a semigroup
$\widetilde{L}_{ev}\left( M\right) $.
Triples of the form $\left( 1_N,1_N,1_N\right) $ are called
{\em trivial\/}. Two triples $u_1$ and $u_2$ are called {\em equivalent\/}
if for some trivial triple
$u_0$ there exists a homotopy
\[
u_1\oplus u_0\sim u_2\oplus u_0.
\]
We consider the Grothendieck group generated by the abelian semigroup of
classes of equivalent triples:
\[
L_{ev}\left( M\right) \stackrel{{\rm def}}{=}K\left( \left. \widetilde{L}%
_{ev}\left( M\right) \right/ \approx \right) .
\]
The index of a triple extends to a homomorphism of abelian groups
\[
\limfunc{ind}_a:L_{ev}\left( M\right) \longrightarrow {\bf Z.}
\]
Let us define the following
two functionals on the group $L_{ev}\left( M\right) $:
$$
\begin{array}{lll}
\limfunc{ind}_t\left( D,P_1,P_2\right) & = & \limfunc{ind}%
_t\left( \sigma \left( D\right) :\func{Im}\sigma \left( P_1\right)
\rightarrow \func{Im}\sigma \left( P_2\right) \right) , \\[1mm]
d\left( D,P_1,P_2\right) & = & d\left( P_1,P_2\right).
\end{array}
$$
Consider the mapping
\begin{eqnarray}
& & d:L_{ev}\left( M\right) \longrightarrow {\bf Z}\left[ \frac 12\right],
\nonumber\\
& & d\left( D,P_1,P_2\right)
\stackrel{{\rm def}}=d\left( \left[ P_1\right] -\left[ P_2\right]
\right).\label{zvezda}
\end{eqnarray}
It is defined by the projections
$P_1,P_2$ and is independent of the operator acting between them. Moreover,
the homomorphism (\ref{zvezda}) is independent of the normalization
$\chi :K\left(M\right) \longrightarrow {\bf R}.$
Let us verify the last statement. Without loss of generality, we can
assume that the vector bundles ${\rm Im}\sigma'(P_{1,2})$ are
pullbacks from the manifold $M$:
$$
{\rm Im}\sigma'(P_{1,2})=\pi^*_PE_{1,2}.
$$
Moreover, the vector bundles $E_{1,2}$ are isomorphic (an isomorphism is
given by the elliptic symbol $\sigma(D)$ lowered to $M$ with the help of a
nonsingular vector field on the manifold). Thus, the $K(M)$ components (in
the sense of the sequence (\ref{seq1})) are equal for the two projections
$P_{1,2}$.

Let us now construct a topological invariant of triples that extends
to a homomorphism of groups
\[
\limfunc{ind}_t:L_{ev}\left( M\right) \longrightarrow
{\bf Q}.
\]
This invariant is determined by the principal symbol of the operator alone.
To this end, consider the quotient space
\[
\widetilde{P}M=\left. \left\{ \left. S^{*}M\times \left[ 0,1\right] \right/
\left( x,\pm \xi ,0\right) \right\} \right/ \left( x,\pm \xi ,1\right)
\]
(see Fig. \ref{ff1}).
\fgr{fig1.eps}{ff1}{}

This space is an oriented manifold with the structure of a fiber bundle
over the manifold $M$. The fiber is formed by two odd-dimensional projective
spaces with small disks deleted\footnote{The projective space with a hole
is a multidimensional analog of the M{\"o}bius band, but in our case it is
orientable.}
\[
P_{-,x}=\left. S_x^{*}M\times \left[ 0,\frac 12\right] \right/ \left( x,\pm
\xi ,0\right) ,P_{+,x}=\left. S_x^{*}M\times \left[ \frac 12,1\right]
\right/ \left( x,\pm \xi ,1\right),
\]
glued along their common boundary
$S_x^{*}M\times \left\{ \frac 12\right\} $
(see Fig. \ref{ff2}):
\fgr{fig2.eps}{ff2}{}
\[
\widetilde{P}M=P_{-}M\bigcup\limits _{S^{*}M}P_{+}M
\]
(the resulting fiber $\widetilde{P}_xM$ is an analog of the Klein bottle).
Consider two natural projections in the bundles with unit closed interval
as a fiber (see Fig. 2)
\[
P_{\mp }M\longrightarrow P^{*}M.
\]
By means of these we can extend the vector bundles $\func{Im}\sigma' \left
( P_{1,2}\right) $ from the space $P^{*}M\subset P_{\mp }M$ to the entire
$P_{\mp }M$, respectively. Then the isomorphism $\sigma(D)$ of these
bundles over the cotangent spheres $S^{*}M$ enables us to glue them; in
this way we obtain a vector bundle over $\widetilde{P}M$:
\[
\sigma \left( D,P_1,P_2\right) \in \limfunc{Vect}\left( \widetilde{P}%
M\right) ,\qquad \left[ \sigma \left( D,P_1,P_2\right) \right] \in K\left(
\widetilde{P}M\right) .
\]
It is possible to define an analog of the usual "topological" index
on the group $K\left( \widetilde{P}M\right) $
(cf. \cite{ABP1}):
\begin{equation}
\label{beta}
{\rm ind}_t\left( D,P_1,P_2\right) =\left\langle \limfunc{ch}\left[ \sigma \left(
D,P_1,P_2\right) \right] \limfunc{Td}\left( T^{*}M\otimes {\bf C}%
\right) ,\left[ \widetilde{P}M\right] \right\rangle.
\end{equation}

In this notation we prove the following index formula, which is the
main result of the present paper.

\begin{theorem}
Let $u=\left( D,P_1,P_2\right) $ be an elliptic operator
acting in the subspaces defined by even pseudodifferential projections
on an odd-dimensional manifold $M.$ Then the following index formula is
valid:
\begin{equation}
\limfunc{ind}_au=\limfunc{ind}_tu+d\left( u\right) .
\label{findx}
\end{equation}
\end{theorem}

\noindent {\em Proof\/}. The left- and right-hand sides of~(\ref{findx})
define homomorphisms of groups
\[
L_{ev}\left( M\right) \longrightarrow {\bf Q.}
\]
To prove the coincidence of these homomorphisms, by virtue
of the absence of torsion in the group of rationals, it suffices
to check the equality only on the elements of the group
$L_{ev}\left( M\right) $ that rationally generate it.
The exactness of the sequence (\ref{seq1}) together with the fact that
the group $K\left( P^{*}M\right) /K\left(M\right) $ is purely a
torsion group imply that for an arbitrary triple
$\left( D,P_1,P_2\right) $ the projections $P_{1,2}$
are rationally homotopic (i.e. a direct sum of the form
$P_{1,2}\oplus\cdots\oplus P_{1,2}$ is homotopic)
to projections that differ from projections on the spaces of sections of
vector bundles by finite rank projections. Hence, as the
elements rationally generating the group $L_{ev}\left( M\right)$,
we can take triples of the form
\[
u=\left( \left(
\begin{array}{cc}
0 & 0 \\
D & 0
\end{array}
\right) ,1_E\oplus n, m\oplus 1_F\right) ,
\]
where $E,F\in \limfunc{Vect}\left( M\right) ,$  $n$ and $m$ are
finite rank projections, and the operator
\[
D:C^\infty \left( M,E\right) \longrightarrow C^\infty \left( M,F\right)
\]
is a usual pseudodifferential operator. Let us show
that formula (\ref{findx}) is true for operators of this form.
Consider the relations
\begin{eqnarray*}
\limfunc{ind}_au &=&\limfunc{ind}D+n-m, \\
d\left( u\right)  &=&n-m+d\left( \left[ 1_{C^\infty \left( M,E\right)
}\right] -\left[ 1_{C^\infty \left( M,F\right) }\right] \right) =n-m
\end{eqnarray*}
(in the second one we use a nonsingular vector field to obtain an
isomorphism of the bundles $E$ and $F$). Thus, to check
formula (\ref{findx}) for a triple $u$, it suffices to
verify that
\begin{equation}
\label{trebform}
\limfunc{ind}D=\limfunc{ind}_t\left( D,1_{C^\infty \left(
M,E\right) },1_{C^\infty \left( M,F\right) }\right) .
\end{equation}

Let us show that the right-hand side of this formula is a slight
modification of the usual Atiyah--Singer index formula for the
operator $D$. We rewrite the Atiyah--Singer index formula in
the form (see \cite{ABP1})
\begin{equation}
\label{alpha1}
\limfunc{ind}D=\left\langle \limfunc{ch}\left[ \sigma \left(
D\right) \right]  \pi  ^{*}\limfunc{Td}\left(  T^{*}M\otimes
{\bf C%
}\right) ,\left[ \widetilde{B}M\right] \right\rangle ,\quad \widetilde{B}%
M=B_{-}M\bigcup\limits_{S^{*}M}B_{+}M,
\end{equation}
where $B_{\pm }M$ is the unit cotangent ball bundle for the manifold
$M.$  Next, consider the two expressions in
(\ref{alpha1}) and (\ref{beta}) as integrals over the respective manifolds
of characteristic classes represented via differential forms
by means of connections in the corresponding vector bundles. The
manifolds $\widetilde{B}M$ and $\widetilde{P}M$ are diffeomorphic
in a neighborhood of $S^{*}M\times \left\{ \frac 12\right\} .$ We note
also that these manifolds carry the orientation-reversing involution
\[
\left( x,\xi ,t\right) \longrightarrow \left( x,-\xi ,t\right)
\]
respecting their parts $P_{\pm}$ and $B_{\pm}$. Consequently, for a
connection in a vector bundle $\sigma \left( D,P_1,P_2\right)$ over
$\widetilde{P}M$ that is invariant under the involution outside of a
certain neighborhood of the cotangent spheres $S^{*}M\times
\left\{1/2\right\} $, there is no contribution to the formula
(\ref{trebform}) from the corresponding domains (since the integrands are
invariant under an orientation-reversing involution); a similar
cancellation happens on $\widetilde{P}M$. The remaining contribution,
coming from the integration over a neighborhood of $S^{*}M\times \left\{1/
2\right\},$ is the same for the two formulas, since the integrands
coincide pointwise.

Thus, the index formula is proved for the case of projections with unit
principal symbol. Such operators, as was noted above, rationally generate
the whole group $L_{ev}(M)$. Hence, the index formula, as well as the
theorem, is proved for the general case.

Let us make two important remarks concerning the topological
term $\limfunc{ind}_t$ in the index formula.

\begin{remark}
{\em In the definition of the group of stably homotopy equivalent
triples
}${L}_{ev}\left( M\right) ${\em \
we could have gone further by factorizing this group by the triples
}$\left( P,P,P\right) $%
{\em \ with an arbitrary even projection
 }$P$, {\em which obviously do not contribute to the index. Let us
 denote the resulting group by the same symbol $L_{ev}$. In this case
we must replace the group
}$K\left( \widetilde{P}M\right)$ {\em (where the principal symbol
of the problem lies) by the quotient
}$\left. K\left( \widetilde{P}M\right) \right/
K\left( P^{*}M\right) $,{\em \ and the corresponding
principal symbol mapping becomes}
\[
L_{ev}\left( M\right) \longrightarrow \left. K\left( \widetilde{P}M\right)
\right/ K\left( P^{*}M\right) .
\]
{\em The last group is an analog of the }$K${\em-functor corresponding
to the difference construction of the usual elliptic theory
}$K_c\left( T^{*}M\right)${\em \ in view of the natural isomorphism
}
\begin{equation}
\left. K\left( \widetilde{P}M\right) \right/ K\left( P^{*}M\right) \simeq
K_c\left( ^{\text{ }p}T^{*}M\right) ,\;\mbox{\em where }^{\text{ }%
p}T^{*}M=P_{-}M\bigcup\limits_{S^{*}M}\left\{ T^{*}M\cap \left\{ \left| \xi \right|%
\geq 1\right\} \right\} .  \label{iso1}
\end{equation}
{\em The topological index}
\[
K_c\left( ^{\text{ }p}T^{*}M\right) \longrightarrow {\bf Q}
\]
{\em is given in this case by the same formula (\ref{beta}).}
\end{remark}

This statement in fact follows from the isomorphism
(\ref{iso1}), which, in turn, is geometrically obvious:
the quotient group on the left-hand side is isomorphic to
the relative group
\[
\left. K\left( \widetilde{P}M\right) \right/ K\left( P^{*}M\right) \simeq
K\left( \widetilde{P}M,P_{+}M\right) ,
\]
and the noncompact spaces  $\widetilde{P}M\setminus P_{+}M$ and
$^{\text{ }p}T^{*}M$ are properly homeomorphic
(i.e., there exists a homeomorphism given by a proper map).

\begin{remark}
{\em \ The topological index }$\limfunc{ind}_t${\em \ of a triple
can be reduced to the topological index of a usual elliptic operator.
Namely, for a triple
}$\left( D,P_1,P_2\right) ${\em \
 }
{\em consider the symbol}
\[
\sigma \left( D\right) ^{-1}\left( x,-\xi \right) \sigma \left( D\right)
\left( x,\xi \right) \oplus 1:\pi _S^{*}E_1\longrightarrow \pi _S^{*}E_1,
\quad \pi_S:S^*M\rightarrow M,
\]
{\em where the direct sum of symbols is taken with respect to the
bundle decomposition}
\[
\pi _S^{*}E=\func{Im}\sigma \left( P_1\right) \left( x,\xi \right) \oplus
\func{Im}\left( 1-\sigma \left( P_1\right) \left( x,\xi \right) \right) .
\]
{\em For this symbol we have}
\begin{equation}
\limfunc{ind}_t\left( D,P_1,P_2\right) =\frac 12\limfunc{ind}%
_t\left( \sigma \left( D\right) ^{-1}\left( x,-\xi \right) \sigma
\left( D\right) \left( x,\xi \right) \oplus 1\right),   \label{halff}
\end{equation}
{\em where the right-hand side is
the topological index of a usual
elliptic operator in spaces.}
\end{remark}

To prove (\ref{halff}), we note that
at the beginning of this section we could have taken
this formula as a definition of
$\limfunc{ind}_t\left( D,P_1,P_2\right)$. Furthermore, in the proof
of the index theorem we would have to check that for a
``classical''\footnote{That is, acting in spaces of sections on the base.}
operator $D$ the expression on the right-hand side in
(\ref{halff}) is equal to the (topological) index of $D$.
Indeed, in this case we obtain
\[
\begin{array}{c}
\limfunc{ind}_t\left( \sigma \left( D\right) ^{-1}\left( x,-\xi
\right) \sigma \left( D\right) \left( x,\xi \right) \oplus 1\right) =%
\limfunc{ind}_t\left( \sigma \left( D\right) ^{-1}\left( x,-\xi
\right) \right) +\limfunc{ind}_t\left( \sigma \left( D\right)
\left( x,\xi \right) \right) = \\
=-\limfunc{ind}_t\left( \sigma \left( D\right) \left( x,-\xi
\right) \right) +\limfunc{ind}_t\left( \sigma \left( D\right)
\left( x,\xi \right) \right) =2\limfunc{ind}_t\left( \sigma \left(
D\right) \left( x,\xi \right) \right) =2\limfunc{ind}D
\end{array}
\]
(in this chain of equalities only the symbols in the classical sense
appear, and ${\rm ind}_t$ is the usual topological index; in the 
next to the last equality we use the orientation-reversing involution $\left( x,\xi \right)
\rightarrow \left( x,-\xi \right) $ on the manifold $S^{*}M$, and the last
equality is the Atiyah--Singer theorem).

\section{General boundary value problems\label{genbvp}}

Let us briefly recall the definition of general boundary value
problems \cite{ScSS18}.

Let $M$  be a smooth compact manifold with boundary. Consider an
$m$th-order elliptic differential operator
\[
D:C^\infty \left( M,E\right) \rightarrow C^\infty \left( M,F\right)
\]
acting in the spaces of sections of vector bundles over $M$.
For some collar neighborhood of the boundary with normal coordinate
$t$ (the interior of the manifold corresponds to positive values of $t$),
consider the jets of sections of the bundle $E$ in the normal direction:
\[
j_{\partial M}^{m-1}u=\left( u(0) ,\left. \frac \partial
{\partial t}u\right| _{t=0},\ldots ,\left. \frac{\partial ^{m-1}}{\partial%
t^{m-1}}u\right| _{t=0}\right) :C^\infty \left( M,E\right) \rightarrow
C^\infty \left( \partial M,E^m\right);
\]
here
$$
E^m=\underbrace{E\oplus E\oplus\ldots\oplus E}_{m \quad\mbox{times}}.
$$
A\emph{ general boundary value problem } is a system of
equations of the form
\begin{equation}
\left\{
\begin{array}{lll}
Du=f, &  & \;u\in C^\infty \left( M,E\right) ,\:
f\in C^\infty \left( M,F\right),\\
Bj_{\partial M}^{m-1}u=g\in {\mathrm{}Im}P, &  & \;{\mathrm{}Im}P\subset C^\infty \left(
\partial M,G\right) .
\end{array}
\right.  \label{bvp1}
\end{equation}
Here the pseudodifferential operator $P$
\[
P:C^\infty \left( \partial M,G\right) \rightarrow C^\infty \left( \partial
M,G\right)
\]
is a zero-order projection, $P^2=P$, and the boundary operator $B$
\[
B:C^\infty \left( \partial M,E^m\right) \rightarrow C^\infty \left( \partial
M,G\right)
\]
is also a pseudodifferential operator whose range is contained in the
subspace ${\mathrm{}Im}P:$ $PB=B.$ The boundary value problem (\ref{bvp1})
is denoted by $\left( D,B,P\right) .$ The classical boundary value
problems correspond to the special case of the unit projection $P=1$ in
formula (\ref{bvp1}).

To state the Fredholm criterion for these boundary value
problems, for each point
$\left( x,\xi ^{\prime }\right) $ of the cosphere bundle
$S^{*}\partial M$ of the boundary
we consider the following ordinary differential equation with constant
coefficients on the line ${\mathbf{}R}\ni t$:
\[
\sigma \left( D\right) \left( x,0,-i\frac d{dt},\xi ^{\prime }\right)
u\left( t\right) =0
\]
(here $\sigma \left( D\right) $ is the principal symbol of $D$). Let us
denote the subspace of Cauchy data of solutions that are bounded as
$t\rightarrow \pm \infty $ by $L^{\mp }(x,\xi^{\prime }).$ Then for an
elliptic differential operator $D$ one has the decomposition
\begin{equation}
E_x^m=L^{+}\left( x,\xi ^{\prime }\right) \bigoplus L^{-}\left( x,\xi
^{\prime }\right) .  \label{deco3}
\end{equation}
Moreover, in this case the families $L^{\pm }\left( x,\xi
^{\prime }\right) $ define smooth vector subbundles
\[
L^{\pm }\subset \pi ^{*}E^m,\qquad \mbox{where }\quad\pi :S^{*}\partial
M\rightarrow \partial M.
\]

\begin{remark}
\label{elluse}\emph{It can be shown that the decomposition
(\ref{deco3}) is a necessary and sufficient condition for the
ellipticity of the operator} $D$\emph{\ on the boundary $\partial M$.}
\end{remark}

The following criterion for the Fredholm property
of general boundary value problems \cite{ScSS18} holds.

\begin{theorem}
\label{fred2}
The boundary value problem $\left( D,B,P\right) $
defines a Fredholm operator
\[
\left(
\begin{array}{c}
D \\
Bj_{\partial M}^{m-1}
\end{array}
\right) :H^s\left( M,E\right) \rightarrow
\begin{array}{c}
H^{s-m}\left( M,F\right)  \\
\bigoplus  \\
{\mathrm{}Im}P\subset H^\delta \left( \partial M,G\right)
\end{array}
\]
for $s>m-1/2$
if and only if the operator $D$ is elliptic and the restriction
of the principal symbol of the boundary operator $B$ to the subbundle
$L^{-}$ defines an isomorphism
\[
\sigma \left( B\right) :L^{-}\left( D\right) \rightarrow {\mathrm{}Im}\,
\sigma \left( P\right) .
\]
It is supposed here that the orders of the components of the operator
$B$ with respect to the jet $j_{\partial M}^{m-1}$ are consistent
with the indices $s,\delta$ of Sobolev spaces.
\end{theorem}

The class of boundary value problems for differential
operators is too narrow
for making homotopies of elliptic symbols. It turns out that the following
simple generalization of this class enables us to carry out the necessary
homotopies; hence, we can apply topological methods to the index problem
for the general boundary value problems.

In this paper we will use operators which are differential with respect
to the normal variable $t$. Namely, let us consider operators $D$ on the
manifold $M$ that have the form
\begin{equation}
D=\sum_{k=0}^mA_k\left( t\right) \left( -i\frac \partial {\partial t}\right)
^k
\label{alpha}
\end{equation}
in a collar neighborhood of the boundary,
where $A_k\left( t\right) $ are smooth families of pseudodifferential
operators of orders $m-k,$ the operator
$A_m\left( t\right) $ being a homomorphism of vector bundles.
Second, on the entire manifold the operator $D$ must be represented by
the construction of a pseudodifferential operator with a continuous symbol
(see \cite{AtSi1},\cite{ReSc1} or \cite{Hor3}). Let us note that the necessity to
consider operators with continuous symbols stems from the fact that
the symbol of the operator (\ref{alpha}) already is not smooth in
general.\footnote{ This is the well-known problem of nonpseudodifferentiality
of tensor products (see \cite{AtSi1}). Here it is, for example, the
composition of operators $\frac \partial {\partial t}$ and $%
A_1\left( t\right) .$}

For the operators of the class just defined, the general boundary value
problems are posed in just the same way as in the usual classical case
(\ref{bvp1}). The definitions of the subspaces $L^{\pm }$ and the
criterion for the Fredholm property remain valid (Theorem \ref{fred2}).

Let us consider examples of boundary value problems
for operators of the form (\ref{alpha}). These special boundary
value problems will enable us (in the classical case)
to reduce the index problem to the known case of a closed manifold.

\begin{example}
\emph{Let a vector bundle }$E$\emph{\ in a neighborhood of the boundary
}$\partial M$\emph{\ be decomposed into a sum of two subbundles}
\begin{equation}
\left. E\right| _{U_{\partial M}}=E_{+}\oplus E_{-}.  \label{deco2}
\end{equation}
\emph{For the two bundles}$\left. E_{\pm
}\right| _{\partial M}$\emph{\ consider first-order elliptic
operators }$%
\Lambda _{\pm }$\emph{\
with the principal symbols equal to}
$\left| \xi ^{\prime }\right| .$\emph{\ We also choose a first-order
operator }$\Lambda $%
\emph{\ with principal symbol }$\left| \xi \right|$\emph{
acting in the bundle }$E$\emph{\ on the entire manifold $M$. Let us define,
in conformity with the decomposition (\ref{deco2}), the following
elliptic first order operator in a collar neighborhood of the boundary:
}
\begin{equation}
\label{mu}
D=\left( -i\frac \partial {\partial t}+i\Lambda _{+}\right) \oplus \left(
+i\frac \partial {\partial t}+i\Lambda _{-}\right) :C^\infty \left(
U_{\partial M},E\right) \rightarrow C^\infty \left( U_{\partial M},E\right) .
\end{equation}
\emph{For the principal symbol of this operator on the boundary
we obtain the equality}
\[
L^{-}\left( D\right) =0\oplus \pi ^{*}E_{-}.
\]
\emph{Hence, the following boundary condition for this operator
is elliptic:}
\begin{equation}
\left. u_{-}\right|_{\partial M}=g\in C^{\infty}(\partial M,E_-)
,\qquad\mbox{\em for }u=\left(
u_{+},u_{-}\right) \in C^\infty \left( U_{\partial M},E_{+}\oplus
E_{-}\right) .  \label{bndcond}
\end{equation}
\emph{ This operator can be extended to the entire manifold. To this end,
consider a smooth cutoff function $\chi$ on
}$M$, $0\leq \chi \left( t\right) \leq 1$,\emph{
which is identically equal to }%
$1$\emph{\ for }$0\leq t\leq 1/3$\emph{\ and vanishes for }
$t\geq 2/3.$\emph{\ The desired extension is now given, say,
by the formula}
\begin{equation}
D=\chi \left( t\right) \left[ \left( -i\frac \partial {\partial t}+i\Lambda%
_{+}\right) \oplus \left( +i\frac \partial {\partial t}+i\Lambda _{-}\right)
\right] +\left( 1-\chi \left( t\right) \right) i\Lambda .  \label{primer}%
\end{equation}
\emph{The operator\ }$D$\emph{\ with the boundary condition (\ref
{bndcond}) defines an elliptic boundary value problem. It is
well known that this boundary value problem has index zero.
This fact can be proved by noting that the family of boundary value
problems}
\[
D+ip
\]
\emph{is an elliptic family with parameter
}$p$\emph{\ in the half-plane
${\rm Re}p>0$ in the sense of Agranovich--Vishik \cite{AgVi1}.
}\emph{\ Consequently, it is invertible
for large values of
}$p.$\emph{\ Besides, the invertibility of the family}
$D+ip$\emph{\ can be shown directly (see \cite{Hor3}).
Correcting the operator $D$ by a finite-dimensional operator, we can
suppose that $D$ is invertible itself.}
\end{example}

In the next section we carry out reductions of boundary value problems,
and there it will be more convenient not to consider an explicit
homotopy of the operator $D$ on the manifold $M$, but rather
consider a homotopy of its
restriction to a small neighborhood of the boundary. More precisely, we
start from homotopies of the form
\[
\sum_{k=0}^mA_k\left( t,\tau \right) \left( -i\frac \partial {\partial
t}\right) ^k,
\]
defined for small values of the parameter $t,$ for example, for $t<1.$
Here $\tau \in \left[ 0,1\right]$ is the parameter of the homotopy. It is
easy to note that for such a homotopy in a neighborhood of the boundary,
one can construct a homotopy of the operator $D$ on the entire manifold.
The required homotopy of the operator $D$ for $t\geq 1$ is constant, while
in the remaining part of the collar neighborhood of the boundary it is
equal to (for the cutoff function $\chi \left( t\right) $ from the example
discussed above)
\[
D_\tau =\sum_{k=0}^mA_k\left( t,\tau \chi \left( t\right) \right) \left(%
-i\frac \partial {\partial t}\right) ^k.
\]

\section{From general to spectral boundary value problems\label{genspec}}

In the present section we show that the methods of index theory of
classical boundary value problems \cite{Hor3} enable us to reduce a
general boundary value problem to a certain spectral boundary value
problem in a canonical way. Moreover, the reduction process does not
affect the space of boundary data defined by the projection $P$. We divide
the reduction procedure into several stages.

Let us briefly comment on the corresponding constructions. Steps 1 and 2
of the reduction are auxiliary in the sense that here the boundary
operator $B$ does not change. On the third step, which is the basis of the
construction, we produce a homotopy of the boundary operator $B$ to the
trivial one.

\subsection*{Step 1. Reduction to first-order operators}

The index problem for the boundary value problem
$\left( D,B,P\right) $ for an operator of order
$m$ is reduced to a similar problem for a first-order
operator by the following theorem.

\begin{theorem}
An elliptic boundary value problem $\left( D,B,P\right) $,
$\deg D=m$, is stably homotopic to the boundary
value problem
\[
\left( D^{\prime },B^{\prime },P\right) \circ mD_{+}^{m-1},\qquad \mbox{ }
mD_{+}^{m-1}\stackrel{{\rm def}}{=}{
\underbrace{D_{+}^{m-1}\oplus \ldots \oplus
D_{+}^{m-1}}_{m\mbox{{\footnotesize \em times}}}
},
\]
where $D_{+}$ is the invertible first-order operator from Example
{\em 1} (for the vector bundles $E_{-}=0,E=E_{+}$) that does not
require boundary conditions and $
\left( D^{\prime },B^{\prime },P\right) $ is a boundary value
problem for a first order operator. In addition, the projection
$P$ is constant in the homotopy.
\end{theorem}

\noindent {\em Proof\/}.
Consider the direct sum of boundary value problems
\begin{equation}
\left( D,B,P\right) \oplus \bigoplus\limits_1^{m-1}D_{+}^m.  \label{aux1}
\end{equation}
Its index coincides with the index of the original problem,
since the operator $D_{+}$ is invertible. Let us represent
the operator $D$ in the form
\[
D=\sum\limits_{k=0}^mD_k\left( t\right) \left( -i\frac \partial {\partial
t}-i\Lambda \right) ^k\left( -i\frac \partial {\partial t}+i\Lambda \right)
^{m-k}\equiv \sum\limits_{k=0}^mD_k\left( t\right) D_{-}^kD_{+}^{m-k},
\]
where $\Lambda $ is an operator with principal symbol
$\left| \xi ^{\prime }\right| $ acting in the vector bundle
$E$ in a neighborhood of the boundary. Consider the following homotopy
of operators (in these matrices, we omit the argument
$t$ of the coefficients $D_k\left( t\right)$ for brevity)
\begin{eqnarray*}
D_\tau &=&\left(
\begin{array}{cccc}
D-\tau ^mD & 0 & \ldots & 0 \\
0 & D_{+}^m & \ldots & 0 \\
0 & 0 & D_{+}^m & \ldots \\
0 & 0 & \ldots & D_{+}^m
\end{array}
\right) + \\
&&\left(
\begin{array}{cccc}
\tau ^mD_0D_{+}^m & \tau ^{m-1}D_1D_{+}^m & \ldots & \tau
D_{m-1}D_{+}^m+\tau D_mD_{-}D_{+}^{m-1} \\
-\tau D_{-}D_{+}^{m-1} & 0 & \ldots & 0 \\
0 & -\tau D_{-}D_{+}^{m-1} & \ldots & \ldots \\
0 & 0 & \ldots & 0
\end{array}
\right).
\end{eqnarray*}
At the initial point $\tau =0$ we have
$D_{\tau =0}=D\oplus \bigoplus\limits_1^{m-1}D_{+}^m.$
On the other hand, for $\tau =1$ the factorization
required in the theorem is obtained:
\[
D_{\tau =1}=\left(
\begin{array}{cccc}
D_0\left( t\right) D_{+} & D_1\left( t\right) D_{+} & \ldots &
D_{m-1}\left( t\right) D_{+}+D_m\left( t\right) D_{-} \\
-D_{-} & D_{+} & 0 & 0 \\
0 & - D_{-} & D_{+} & \ldots \\
0 & 0 & \ldots & D_{+}
\end{array}
\right) \circ mD_{+}^{m-1}. %
\]
Let us verify that the operator $D_\tau$ is elliptic for
$\tau \in \left[ 0,1\right] $. First of all, we calculate
the subspace $L^{-}\left( D_\tau
\right) .$ Let $U=\left( U_0\left( t\right) ,U_1\left(
t\right) ,\ldots ,U_{m-1}\left( t\right) \right) $ be a bounded
solution as $%
t\rightarrow +\infty $ of the equation
\begin{equation}
\label{eqq}
\sigma \left( D_\tau \right) \left( x,0,-i\frac d{dt},\xi ^{\prime }\right)
U=0.
\end{equation}
We would like to note that the derivatives of the solution $U$ are
bounded too, since it is a solution of an ordinary differential equation
(\ref{eqq})
with constant coefficients. Writing out this equation componentwise,
we obtain the
system of equations
\begin{equation}
\left\{
\begin{array}{l}
\left( 1-\tau ^m\right) \sigma \left( D\right) \left( -i\frac d{dt}\right)
U_0+\left( -i\frac d{dt}+i\right) ^m\left( d_0\tau ^mU_0+\ldots +d_{m-1}\tau
U_{m-1}\right) + \\
+\tau \left( -i\frac d{dt}+i\right) ^{m-1}\left( -i\frac d{dt}-i\right)
d_mU_{m-1}=0, \\
\left( -i\frac d{dt}+i\right) ^mU_j=\tau \left( -i\frac d{dt}+i\right)
^{m-1}\left( -i\frac d{dt}-i\right) U_{j-1}\qquad \mbox{for }%
0<j<m
\end{array}
\right.  \label{krut}
\end{equation}
(here by $d_j$ we denote the principal symbols of the operators
$D_j(0)$). Since the equation
\[
\left( -i\frac d{dt}+i\right) u=0
\]
has no bounded solutions on the half-line, we can cancel this operator in the system
(\ref{krut}) in all equations except the first one. Hence, we obtain
\[
\left( -i\frac d{dt}+i\right) U_j=\tau \left( -i\frac d{dt}-i\right)
U_{j-1}.
\]
Successively substituting these relations into one another, we obtain
\[
\left( -i\frac d{dt}+i\right) ^jU_j=\tau ^j\left( -i\frac d{dt}-i\right)
^jU_0.
\]
Taking into account these equations in the first equation in
(\ref{krut}), we find that it is reduced to the requirement
\begin{equation}
\sigma \left( D\right) \left( -i\frac d{dt}\right) U_0=0.  \label{reduk}
\end{equation}
This implies the ellipticity of the operator
$D_\tau :$ equation (\ref{reduk}) on the entire line
has no bounded solutions, i.e. the decomposition
(\ref{deco3}) is valid.

Thus, the following description of the subbundle
$L^{-}\left( D_\tau \right)$ is obtained. The projection on the first
term in the sum
\[
E\oplus \bigoplus\limits_1^{m-1}E\stackrel{pr}{\longrightarrow }E
\]
induces an isomorphism of vector bundles
\[
L^{-}\left( D_\tau \right) \stackrel{pr}{\longrightarrow }L^{-}\left(
D\right) ;
\]
the preimage of an element $u\in L^{-}\left( D\right) $ under
this map is given by the formula
\begin{eqnarray}
U &=&\left( U_0,\ldots ,U_{m-1}\right),  \label{tozhd} \\
U_0 &=&u,  \nonumber \\
\left( -i\frac d{dt}+i\right) ^jU_j &=&\tau ^j\left( -i\frac d{dt}-i\right)
^jU_0.  \nonumber
\end{eqnarray}
Let us decompose the boundary operator in the same way
as the operator $D$ has been rewritten in the above:
\[
Bj_{\partial M}^{m-1}=\left. \sum\limits_{k=0}^{m-1}B_k\left( -i\frac
\partial {\partial t}-i\Lambda \right) ^k\left( -i\frac \partial {\partial
t}+i\Lambda \right) ^{m-1-k}\right| _{t=0}.
\]
Therefore, the principal symbol of the boundary condition
\[
Bj_{\partial M}^{m-1}\circ pr:C^\infty \left( M,E\oplus%
\bigoplus\limits_1^{m-1}E\right) \rightarrow C^\infty \left( \partial
M,G\right)
\]
on the subspace $L^{-}\left( D_{\tau =1}\right) $ has the
factorization (by virtue of (\ref{tozhd}))
\[
\sigma \left( B^{\prime }\right) j_{\partial M}\circ \left( -i\frac%
d{dt}+i\right) ^{m-1},\qquad \mbox{where\quad }\left( \sigma
\left( B^{\prime }\right) j_{\partial M}\right)
U=\sum\limits_{k=0}^{m-1}b_kU_k\left( 0\right) .
\]
Hence, the homotopy
\[
\left( D_\tau ,B\circ pr,P\right) ,\qquad \tau \in \left[ 0,1\right],
\]
of elliptic boundary value problems
connects the initial problem (\ref{aux1}) with the boundary value problem
that is equal to the composition indicated in the theorem:
\[
\left( D_{\tau =1},B\circ pr,P\right) =\left( D^{\prime },B^{\prime%
},P\right) \circ mD_{+}^{m-1}. %
\]
This completes the proof of the theorem.

\subsection*{Step 2. An expression for the operator on the boundary
via the Calder{\'o}n projection}

At this step we show that a first order operator is reduced in a collar
neighborhood of the boundary to an operator that is uniquely determined by
the subbundle $L^{-}\left( D\right).$

Consider a boundary value problem for the first-order elliptic operator
\[
D=-iA_1\frac \partial {\partial t}+A_0:C^\infty \left( M,E\right)
\rightarrow C^\infty \left( M,F\right).
\]
The coefficient $A_1$ is by assumption an isomorphism of bundles
in a neighborhood of the boundary. With the help of this isomorphism,
we identify the bundles $E$ and $F.$ Then the operator $D$ becomes
\[
D=-i\frac \partial {\partial t}+A.
\]
The ellipticity of $D$ implies that the eigenvalues of the principal
symbol of the operator $A$ for nonzero values of the cotangent variable
$\xi ^{\prime }$ on the boundary have nonzero imaginary parts. A simple
calculation shows that the subbundle $L^{-}\left( D\right) $ is generated
by the spectral subspaces of the symbol $\sigma \left( A\right) $ with
negative imaginary parts of the eigenvalues. Let us denote the
corresponding spectral projection by $q$:
\[
q:\pi ^{*}E\rightarrow \pi ^{*}E,\qquad {\mathrm{}Im}q=L^{-}\left( D\right) .
\]
Let us also consider an arbitrary pseudodifferential operator
$Q$ (not necessarily a projection) with principal symbol
$q.$ Such operators are called
\emph{Calder{\'o}n projections\/} for the operator  $D$ (see \cite{Hor3}).

Consider a linear homotopy of operators
\[
D_\tau =-i\frac \partial {\partial t}+\left( 1-\tau \right) A-i\tau \left(
2Q-1\right) \Lambda ,
\]
where $\Lambda $  is an operator with principal symbol
equal to $\left| \xi ^{\prime }\right| .$ Since the projection
$q\left( \xi ^{\prime }\right) $ is the spectral projection for the
symbol $\sigma \left( A\right)
\left( \xi ^{\prime }\right) ,$ we see
that the eigenvalues of the symbol
\[
\left( 1-\tau \right) \sigma \left( A\right) \left( \xi ^{\prime }\right)
-i\tau \left( 2q-1\right) \left| \xi ^{\prime }\right|
\]
are, respectively, equal to
(here $\lambda $  is an arbitrary eigenvalue of the symbol
$\sigma \left( A\right) $)
\[
\left( 1-\tau \right) \lambda +i\tau {\mathrm sign}\left(
{\mathrm{}Im}\lambda \right)  .
\]
This implies the ellipticity of the operator $D_\tau $ and the
independence of the subbundle $L^{-}\!\!\left( D_\tau \right)$ of the
parameter $\tau .$

In this way, we present a homotopy of boundary value problems
\[
\left( D_\tau ,B,P\right),\quad \tau\in [0,1].
\]
As a result, the principal symbol of the operator $D$ obtained
is determined on the boundary by the principal symbol of the
Calder{\'o}n projection $Q$
\begin{equation}
D_{\tau =1}=-i\frac \partial {\partial t}-i\left( 2Q-1\right) \Lambda .
\label{after2}
\end{equation}
At the end of this step, let us consider the orthogonal projection
$q^{\prime }$ onto the subbundle $L^{-}\left( D\right) $ and
a pseudodifferential projection $Q^{\prime }$ with principal symbol
$q^{\prime }.$ Then a linear homotopy
\[
-i\frac \partial {\partial t}-i\left( 2\left( \tau Q^{\prime }+\left( 1-\tau
\right) Q\right) -1\right) \Lambda
\]
leads for $\tau =1$ to an operator
\[
-i\frac \partial {\partial t}-i\left( 2Q^{\prime }-1\right) \Lambda
\]
with principal symbol uniquely determined by the subbundle
$L^{-}\left( D\right).$ Hence, it can be assumed in what follows that
the symbol $q$ of the Calder{\'o}n projection is the orthogonal projection
onto the subbundle $L^-(D)$.

\subsection*{Step 3. Reduction of the Calder{\'o}n projection
to the projection of the boundary data}

In this section we make a homotopy of the operators $D$ and $B$ of the
boundary value problem $\left(D,B,P\right)$. As a result of this homotopy,
the boundary operator $B$ is transformed from an operator with principal
symbol giving an isomorphism of subbundles
\[
L^{-}\left( D\right) \stackrel{\sigma \left( B\right) }{\longrightarrow }%
{\mathrm{}Im}\,\sigma \left( P\right)
\]
to the identity operator. The bundle $L^{-}\left( D\right) $, in
particular, is deformed into ${\mathrm{}Im}\,\sigma \left(P\right) .$
Formula (\ref{after2}) shows that to construct such a homotopy of boundary
value problems, it suffices to produce a homotopy of the principal symbol
of the Calder{\'o}n projection $q$, as well as of the principal symbol of
the boundary operator $\sigma \left( B\right) .$ From the theory of vector
bundles it is known that isomorphic subbundles (in our case these are
$L^-(D)$ and ${\rm Im}\sigma(P)$) can be deformed into one another as
subbundles in the ambient vector bundle (provided the dimension of this
bundle is large enough). In this particular situation, let us write down
an explicit formula for such a homotopy.

Let us realize the bundle $G\in {\mathrm Vect}\left( \partial
M\right) $ as a subbundle of the trivial ${\mathbf{}{C}}^N$
\[
G\stackrel{\alpha }{\subset }{\mathbf{}{C}}^N.
\]
Recall that the boundary data of the boundary value problem lie
in the space of sections of this bundle: ${\mathrm{}Im}P\subset C^\infty
\left( \partial M,G\right) $.

Consider the direct sum of boundary value problems
\[
\left( D,B,P\right) \oplus D_{+},
\]
where an invertible operator $D_{+}$ is the operator from
Example 1 corresponding to the bundles
$E_{+}=E={\mathbf{}{C}}^N.$ The desired homotopy of subbundles
$$
L^{-}\left( D\right) \oplus
0,0\oplus {\mathrm{}Im}\sigma \left( P\right) \subset \pi ^{*}\left( E\oplus
{\mathbf{}{C}}^N\right)
$$
in this notation is  a homotopy of rotation by $90^{\circ }$ with the help
of the symbol $\sigma \left( B\right) .$ More precisely, for an angle
$\tau \in \left[ 0,\pi /2\right] $ we define the subbundle
$$
L^{-}\left( \tau \right)
\subset \pi ^{*}\left( E\oplus {\mathbf{}{C}}^N\right)
$$
whose fiber over a point
$\left( x,\xi^{\prime }\right) \in S^{*}\partial M$ is generated
by vectors of the form
\begin{equation}
v_\tau =\left(\cos \tau v,\sin \tau \sigma \left( B\right) v\right) \qquad %
\mbox{for all }v\in L^{-}\left( D\right) .
\label{grandiozo}
\end{equation}
The homotopy of the principal symbol of the boundary operator $B$
is carried out by the formula
\[
\sigma \left( B_\tau \right) v_\tau \stackrel{{\rm def}}{=}\sigma \left( B\right)
v.
\]
It is clear from (\ref{grandiozo}) that at the end
of the homotopy for $\tau=\pi/2$ we obtain
\[
L^{-}\left( \pi /2\right) =0\oplus {\mathrm{}Im}\sigma \left( P\right) ,\quad
\sigma \left( B_{\pi /2}\right) ={\mathrm Id}.
\]
Thus, the operator $\left( D_{\pi /2},B_{\pi
/2},P\right) $ defines the desired spectral boundary value problem.

This completes the reduction of a general boundary value problem to the
corresponding spectral boundary value problem of the form
(\ref{bvpspec}).

\section{An index formula for even boundary value problems\label{indevn}}

In this section we obtain an index formula for spectral boundary value
problems of the form  (\ref{bvpspec}) under the additional assumption
that the projection $P$ is even.\footnote{In combination with the
reductions of the previous section, this implies
an index formula for boundary value problems of the general form (\ref{bvp})
under the parity condition.}
We prove the index formula by a reduction
to a classical boundary value problem. To this end, let us study the
relationship between classical and even boundary value problems.

Consider the Grothendieck group $K\left( P_{ev}\left( X\right) \right)$
generated by the abelian semigroup of classes of equivalent even
projections (see formula (\ref{omeg})). It has the subgroup generated by
unit projections. The corresponding quotient group is denoted by $\left. K
\left( P_{ev}\left(X\right) \right) \right/ K\left( X\right) .$

A spectral boundary value problem of the form (\ref{bvpspec}) on the
manifold $M$ is denoted by ${\mathcal{}D}=\left( D,P\right).$ Such a
boundary value problem is called \emph{even\/} if the projection $P$ is an
even projection on the boundary $\partial M.$ Let us introduce an
equivalence relation for even boundary value problems. Namely, two
boundary value problems ${\mathcal{}D}_1$ and ${\mathcal{}D}_2$ are called
equivalent if there exists an even homotopy
\[
{\mathcal{}D}_1\oplus D_{+}\sim D_{+}\oplus
{\mathcal{}D}_2,
\]
where $D_{+}$ is, as in the previous section, an invertible operator
from the Example 1 for the choice of bundles $E=E_{+}.$
The Grothendieck group of even boundary value problems
is denoted by $K\left( {\mathcal{}D}_{ev}\left(M\right) \right) .$
It has the subgroup generated by the classical boundary value
problems from Example 1 for arbitrary vector bundles $E$ on $M$
and $E_{+}$ on $\partial M$, respectively. The quotient group is
denoted by
\begin{equation}
\left. K\left( {\mathcal{}D}_{ev}\left( M\right) \right) \right/ \left(
K\left( \partial M\right) \oplus K\left( M\right) \right) .  \label{alph}
\end{equation}
Finally, we need the group
\[
\left. K\left( {\mathcal{}D}\left( M\right) \right) \right/ \left( K\left(
\partial M\right) \oplus K\left( M\right) \right)
\]
generated by classical boundary value problems
$(P=1)$. It is obtained by the same construction as before,
with the replacement of the condition that the projection $P$ is
even by the condition $P=1$. The resulting quotient group, which is
similar to (\ref{alph}), does not require a new notation,
since it can be identified with the usual
$K$-group of vector bundles with compact support
\begin{equation}
\left. K\left( {\mathcal{}D}\left( M\right) \right) \right/ \left( K\left(
\partial M\right) \oplus K\left( M\right) \right) \stackrel{\gamma }{%
\rightarrow }K_c\left( T^{*}\left( M\backslash \partial M\right) \right)
\equiv K_c\left( T^{*}\stackrel{\circ }{M}\right) .  \label{kcomp}%
\end{equation}
Let us recall the definition of the map $\gamma $.
First, we use the isomorphisms
\[
K_c\left( T^{*}\left( M\backslash \partial M\right) \right) \simeq K\left(
B^{*}M,\partial \left( B^{*}M\right) \right) ,\qquad \partial \left(
B^{*}M\right) =S^{*}M\cup \left. B^{*}M\right| _{\partial M},%
\]
where $B^{*}M$ is the unit ball bundle of the manifold
$M,$ and $\partial \left( B^{*}M\right) $ is its full boundary.
This allows us to identify the elements of the group
$K_c\left( T^{*}\!\stackrel{\circ }{M}\right) $ with the
help of the difference construction both with isomorphisms
of vector bundles $E,F\in {\mathrm{}Vect}\left(M\right) $
\begin{equation}
\sigma :\pi ^{*}E\rightarrow \pi ^{*}F,\quad \pi :S^{*}M\rightarrow M,
\label{tmp11}
\end{equation}
that do not depend on $\xi $ in a neighborhood of the boundary
and with isomorphisms defined everywhere on
$\partial \left( B^{*}M\right) .$

Consider a classical boundary value problem $\left( D,1\right) $ of the
form (\ref{bvpspec}). In a neighborhood of the boundary it can be
obviously rewritten as the boundary value problem from Example 1. Then the
homomorphism $\gamma $ is by definition equal to
\begin{equation}
\label{vsego}
\gamma \left( D,1\right) \stackrel{{\rm def}}{=}\sigma \left( D\right) \left( \xi
\right) +i\chi \left( t\right) ,%
\end{equation}
for a cutoff function $\chi \left( t\right) $ equal to $1$ on $\partial
M$, as above. It follows from (\ref{mu}) that the principal symbol
(\ref{vsego}) is invertible on $\partial \left( B^{*}M\right) .$ It is
also not difficult to construct the inverse mapping for $\gamma $:
\begin{equation}
K_c\left( T^{*}\stackrel{\circ }{M}\right) \stackrel{\gamma ^{\prime }}{
\rightarrow }\left. K\left( {\mathcal{}D}\left( M\right) \right) \right/
\left( K\left( \partial M\right) \oplus K\left( M\right) \right) .
\label{iso}
\end{equation}
To an isomorphism $\sigma $ (see (\ref{tmp11})) that is
independent of $\xi $ over a neighborhood of the boundary
$\partial M$, this map assigns a classical boundary value
problem (\ref{bvpspec}) in the following way.
For the symbol $\sigma $ we construct an elliptic first-order
pseudodifferential operator that has the form
\[
D=\sigma \circ \Lambda
\]
near the boundary, where $\sigma =\sigma \left( x\right) $ is a homomorphism
of vector bundles and $\Lambda $ has the principal symbol
$\left| \xi \right| .$ It remains to modify the operator
$D$ near the boundary, as it was done in Example 1:
\[
D^{\prime }=\sigma \circ \left( \left( 1-\chi \left( t\right) \right)%
\Lambda +\chi \left( t\right) \left( -\frac \partial {\partial t}+\Lambda%
_{+}\right) \right).
\]
The operator $\Lambda _{+}$ here has the principal symbol $\left| \xi
^{\prime }\right| .$ Finally, we define
\[
\gamma ^{\prime }\left( \sigma \right) \stackrel{{\rm def}}{=}\left[ D^{\prime
}\right] \in \left. K\left( {\mathcal{}D}\left( M\right) \right) \right/
\left( K\left( \partial M\right) \oplus K\left( M\right) \right).
\]
This is well defined, since the operator $D^{\prime }$ defines a Fredholm
boundary value problem without boundary conditions.

\begin{remark}
\label{indclass}
\emph{In terms of the isomorphism (\ref{iso}) it is easy to prove the
following index formula for classical boundary value problems of the form
(\ref{bvpspec}):}
\[
{\mathrm{}ind}\left( D,1\right) =p_{!}\gamma \left( D,1\right) ,\qquad p:%
\stackrel{\circ }{M}\rightarrow pt,%
\]
{\em
where $p_!:K_c\left(T^*\stackrel{\circ}{M}\right)\rightarrow K(pt)={\bf Z}
$ is the direct image in $K$-theory.
}
\end{remark}

\noindent {\em Proof\/}. Indeed, consider a boundary value problem
$\left(D^{\prime },1\right) $ on the manifold $M$ that coincides with the
original problem $\left( D,1\right) $ in a neighborhood of the boundary
and on the whole manifold is a boundary value problem from Example 1. Let
us examine the composition
\[
D_0=\left( D^{\prime },1\right) ^{-1}\circ \left( D,1\right).%
\]
Its principal symbol $\sigma(D_0)$ is the identity isomorphism
over a neighborhood of the boundary $\partial M$. Thus,
the index of the operator $D_0$ can be computed by the Atiyah--Singer formula
\[
{\mathrm{}ind}D_0=p_{!}\left[ \sigma \left( D_0\right) \right] ,\qquad \left[
\sigma \left( D_0\right) \right] \in K_c\left( T^{*}\stackrel{\circ }{M}
\right) .
\]
The boundary value problem $\left( D^{\prime },1\right) $, however, has
index zero and also defines the trivial element in the $K$-group:
\[
\gamma \left( D^{\prime },1\right)
=0\in K_c\left( T^{*}\stackrel{\circ }{M}
\right) .
\]
Hence, we obtain the desired formula
\[
{\mathrm{}ind}\left( D,1\right) ={\mathrm{}ind}D_0=p_{!}\left[ \sigma \left(
D_0\right) \right] =p_{!}\gamma \left( D,1\right) .
\]

\begin{remark}
\emph{The  reductions of the previous section give an isomorphism
of the group (\ref{kcomp}) and the group of stable homotopy classes
of classical boundary value problems (\ref{bvp}).
In this way, considering spectral boundary value problems
(\ref{bvpspec}), we make no loss of generality and cover the general
case as well.}
\end{remark}

The three groups introduced above are related by an exact sequence.

\begin{proposition}
The sequence
\begin{equation}
K_c\left( T^{*}\stackrel{\circ }{M}\right) \stackrel{\alpha }{\rightarrow }
\left. K\left( {\mathcal{}D}_{ev}\left( M\right) \right) \right/ \left(
K\left( \partial M\right) \oplus K\left( M\right) \right) \stackrel{\beta }{%
\rightarrow }\left. K\left( P_{ev}\left( \partial M\right) \right) \right/
K\left( \partial M\right)   \label{exact}
\end{equation}
is exact. Here $\alpha $ is induced by the embedding
of classical boundary value problems into even ones and
the map $\beta $ is induced by the forgetful map
\[
\left( D,P\right) \longrightarrow P.
\]
\end{proposition}

\noindent {\em Proof\/}. The equality $\beta \circ \alpha =0$ is obvious,
since the projection in the classical boundary value problem is the unit
projection and hence defines the trivial element in the group $\left.
K\left( P_{ev}\left( \partial M\right) \right) \right/ K\left(\partial
M\right) .$

Let us now verify the inclusion $\ker \beta \subset {\mathrm{} Im}\alpha .
$ Suppose that for an even boundary value problem $\left(D,P\right) $ one
has
\[
\beta \left[ D,P\right] =0\in \left. K\left( P_{ev}\left( X\right) \right)
\right/ K\left( X\right) .
\]
This means that there exists a homotopy of even projections that
connects the projection $P$
and the projection on the space of sections
of a vector bundle on $\partial M$, denoted by
$P^{\prime }=P_{C^\infty \left( \partial M,G^{\prime }\right) }$.
Let us denote this homotopy by $P_t:$ $P_0=0,$
$P_1=P^{\prime },$ and lift this homotopy of projections
to a homotopy of spectral boundary value problems
$\left( D_t,P_t\right) .$ In the Grothendieck group
$\left. K\left( {\mathcal{}D}_{ev}\left( M\right)
\right) \right/ \left( K\left( \partial M\right) \oplus
K\left( M\right) \right) $, we have
\[
\left[ D,P\right] =\left[ D_1,P_1\right] =\alpha \left[ D_1,1_{C^\infty
\left( \partial M,G^{\prime }\right) }\right] ,
\]
since the boundary value problem  $\left( D_1,P_1\right) $
is classical.

This establishes the exactness of the sequence (\ref{exact}).

From now on we assume that the manifold $M$ is \emph{even-dimensional\/}.
The third term in the sequence (\ref{exact}) is simplified in this case:
according to (\ref{seq1}), finite-dimensional projections generate a
subgroup in this term that is isomorphic to $\bf Z$
$$
{\bf Z}\subset\left. K\left( P_{ev}\left( X\right) \right) \right/
 K\left( X\right),
$$
and the quotient group
\[
\left. K\left( P_{ev}\left( X\right) \right) \right/ \left( K\left( X\right)
\oplus {\mathbf{}Z}\right) \simeq K\left( P^{*}X\right)/K(X)
,\quad P^{*}X\mbox{ is the projectivization of }S^{*}X,
\]
according to \cite{Gil7}, consists of elements of finite orders that are
powers of 2.

In a similar fashion, in the Grothendieck group of even
boundary value problems there is a subgroup
\[
{\mathbf{}Z\subset }\left. K\left( {\mathcal{}D}_{ev}\left( M\right) \right)
\right/ \left( K\left( \partial M\right) \oplus K\left( M\right) \right),
\]
which is
generated by boundary value problems for zero operators  with
finite rank projections on the right-hand sides.
This enables us to refine the sequence (\ref{exact}):
\begin{equation}
K_c\left( T^{*}\stackrel{\circ }{M}\right) \stackrel{\alpha }{\rightarrow }
\left. K\left( {\mathcal{}D}_{ev}\left( M\right) \right) \right/ \left(
K\left( \partial M\right) \oplus K\left( M\right) \oplus {\mathbf{}Z}\right)
\stackrel{\beta }{\rightarrow }K\left( P^{*}\left( \partial M\right) \right)
/K\left( \partial M\right) .  \label{exact1}
\end{equation}

For an even boundary value problem, we define an analog of the topological
index of Atiyah--Singer. Let us consider the double
\[
2M\stackrel{{\rm def}}{=}M\bigcup\limits_{\partial M}M
\]
of the manifold $M$. On this manifold we have the (continuous) elliptic
symbol
\[
\sigma(\tilde D)(\xi)=
\sigma \left( D\right) \left( \xi \right) \cup \sigma \left( D\right) \left(%
-\xi \right),
\]
on the first copy of the manifold this symbol is equal to
the original symbol $\sigma(D)(\xi)$, and on the second copy it is
$\sigma(D)(-\xi)$. The continuity at the place of gluing
follows from the equality
\[
\sigma \left( D\right) \left( -\tau ,\xi ^{\prime }\right) =\sigma \left(
D\right) \left( -\tau ,-\xi ^{\prime }\right) .
\]
In this notation the following index formula for spectral boundary
value problems is valid.

\begin{theorem}
\label{thmain}
\begin{equation}
{\mathrm{}ind}\left( D,P\right) =\frac 12{\mathrm{}ind}_t\left(
\sigma \left(\tilde D\right)
\right) -d\left( P\right),  \label{main}
\end{equation}
where $d(P)$ is the dimension of the projection $P$
for the trivial  normalization (see Theorem \/{\em\ref{thmmn}}).
\end{theorem}

\noindent {\em Proof\/}.
The index of the boundary value problem and the right-hand side
of (\ref{main}) extend to homomorphisms of
abelian groups, denoted, respectively, by
\[
{\mathrm{}ind},{\mathrm{}ind}^{\prime }:\left. K\left( {\mathcal{}D}%
_{ev}\left( M\right) \right) \right/ \left( K\left( \partial M\right) \oplus
K\left( M\right) \right) \longrightarrow {\mathbf{}Q.}
\]
Let us take advantage of the exact sequence (\ref{exact1}). We first check
that the two homomorphisms coincide for the classical boundary value
problems. In this case we obtain $d\left( P\right) =0,$ since we have the
unit projection and the normalization $\chi$ is taken to be trivial.
Hence, it remains to verify the validity of the index formula for the
classical boundary value problems
\begin{equation}
{\mathrm{}ind}\left( D,1\right) =\frac 12{\mathrm{}ind}_t\left(
\sigma \left( D\right) \left( \xi \right) \cup \sigma \left( D\right) \left(%
-\xi \right) \right) .  \label{question}
\end{equation}
By virtue of the isomorphism (\ref{iso}), the left- and
right-hand
sides of (\ref{question}) are homomorphisms of groups
\[
K_c\left( T^{*}\stackrel{\circ }{M}\right) \longrightarrow {\mathbf{}Q.}%
\]
The right-hand side of (\ref{question}) is decomposed for the classical
boundary value problems into two terms
\[
{\mathrm{}ind}_t\left( \sigma \left( D\right) \left( \xi \right)
\cup \sigma \left( D\right) \left( -\xi \right) \right) ={\mathrm{}ind}
_t\left( \sigma \left( D\right) \left( \xi \right) \right) +%
{\mathrm{}ind}_t\left( \sigma \left( D\right) \left( -\xi \right)
\right),
\]
\[
\sigma(D)(\xi),\sigma(D)(-\xi)\in K_c\left( T^{*}\stackrel{\circ }{M}\right).
\]
Let us show that on an even-dimensional manifold the two terms
in the last formula are equal. Indeed, in the cohomological form
we obtain
\begin{equation}
{\mathrm{}ind}_t\left( \sigma \left( D\right) \left( -\xi \right)
\right) =\left\langle {\mathrm{}ch}\left[ \sigma \left( D\right) \left( -\xi
\right) \right] {\mathrm{}Td}\left( T^{*}M\otimes {\mathbf{}{'}}\right)
,\left[ T^{*}\stackrel{\circ }{M}\right] \right\rangle .  \label{tmp0}
\end{equation}
The involution on the space $T^{*}M$
\[
\xi \longrightarrow -\xi
\]
preserves its orientation. Then from (\ref{tmp0}) we obtain
\[
{\mathrm{}ind}_t\left( \sigma \left( D\right) \left( \xi \right)
\right) ={\mathrm{}ind}_t\left( \sigma \left( D\right) \left( -\xi
\right) \right) .
\]

The equality in (\ref{question}) now follows from the index formula for
the classical boundary value problems (see Remark \ref{indclass}):
\[
{\mathrm{}ind}\left( D,1\right) ={\mathrm{}ind}_t\left( \sigma \left(
D\right) \left( \xi \right) \right) ,
\]

The homomorphisms $ {\mathrm{}ind}$ and ${\mathrm{}ind}^{\prime }$
obviously coincide on the subgroup ${\mathbf{}Z}$ generated by
boundary value problems for operators equal to zero. In this way, the
difference ${\mathrm{}ind}-{\mathrm{}ind}^{\prime }$ descends to a
homomorphism of the quotient group
\[
{\mathrm{}ind}-{\mathrm{}ind}^{\prime }:\left. \left[ \left. K\left(
{\mathcal{}D}_{ev}\left( M\right) \right) \right/ \left( K\left( \partial
M\right) \oplus K\left( M\right) \oplus {\mathbf{}Z}\right) \right] \right/
{\mathrm{}Im}\alpha \rightarrow {\mathbf{}Q.}
\]
The exactness of the sequence (\ref{exact1}) implies ${\mathrm{}Im }\alpha
=\ker \beta .$ Let us now note that the homomorphism $\beta $ takes values
in the torsion group (see the sequence (\ref{exact1})). Thus, the
difference ${\mathrm{}ind}-{\mathrm{}ind}^{\prime}$ cannot be nontrivial,
since it is defined on the group consisting of elements of a finite order.

The index formula is thereby proved.

Operators in subspaces on closed manifolds and boundary value
problems on manifolds with boundary studied
above are related.

Let $D$ be an operator in subspaces, on a closed manifold $M$
\begin{equation}
D:\func{Im}P_1\rightarrow \func{Im}P_2,\qquad \func{Im}P_{1,2}\subset
C^\infty \left( M,E_{1,2}\right) .
\label{subb}
\end{equation}
Let us assign an elliptic boundary value problem to this operator.

On the cylinder $M\times \left[ 0,1\right] $, consider an elliptic
first-order operator with constant coefficients along the cylinder
\begin{equation}
D^{\prime }=\frac \partial {\partial t}+\left( 1-P_1\right) \Lambda \left(
1-P_1\right) -P_1\Lambda P_1:C^\infty \left( M\times \left[ 0,1\right]
,E_1\right) \rightarrow C^\infty \left( M\times \left[ 0,1\right]
,E_1\right) ,  \label{opp}
\end{equation}
where $\Lambda$ is a positive self-adjoint operator of the first
order on $M$ with principal symbol $\left| \xi \right| .$
On the two bases of the cylinder we impose
the following boundary conditions: one boundary condition
spectral, while the other is general,
\begin{equation}
\left\{
\begin{array}{c}
\left( 1-P_1\right) u\left( 0\right) =g\in \func{Im}\left( 1-P_1\right)
\subset C^\infty \left( M,E_1\right),  \\
Du\left( 1\right) =g^{\prime }\in \func{Im}P_2\subset C^\infty \left(
M,E_2\right) .
\end{array}
\right.   \label{condd}
\end{equation}
Let us note that the ellipticity of the boundary conditions (\ref{condd})
is equivalent to the ellipticity of the operator in the subspaces
(\ref{subb}).

\begin{proposition}
The index of the boundary value problem
{\em (\ref{opp}), (\ref{condd})} is equal to the index of the operator
in subspaces {\em (\ref{subb})}.
\end{proposition}

\noindent{\em Proof\/}.
First, it suffices to prove this proposition for
the case in which the projections
$P_{1,2}$ are orthogonal. In this situation the operator
(\ref{opp}) has the form
\[
D^{\prime }=\frac \partial {\partial t}+A,\quad A\text{ is a
self-adjoint operator,}
\]
which makes it possible to reduce the proof of the equality of the indices
to a direct calculation (with the help of the eigenfunctions of $A$).
These calculations are omitted.

\section{Eta invariants and even projections\label{defadm}}

Following \cite{Gil7}, we say that a classical pseudodifferential
operator $A$ of integer positive order  $m$
is {\em admissible\/} if its complete symbol
\[
a\left( x,\xi \right) \sim a_m\left( x,\xi \right) +a_{m-1}\left( x,\xi
\right) +...,
\]
satisfies the parity conditions
\[
a_\alpha \left( x,-\xi \right) =\left( -1\right) ^\alpha a_\alpha \left(
x,\xi \right) ,\qquad \forall \xi \neq 0,\quad x,\quad\alpha=m,m-1,m-2,\ldots
\]
(the admissibility of an operator is independent of the choice of a
coordinate system where the complete symbol is considered). We also recall
that the $\eta$-{\em function\/} of a self-adjoint operator $A$ is defined
as
\[
\eta \left( s,A\right) =\sum_{\lambda \in \limfunc{spec}A}\limfunc{sign}%
\lambda \left| \lambda \right| ^{-s},
\]
where the sum is taken over the nonzero eigenvalues of the operator $A$
(with regard to their multiplicities). It is well known that for an
elliptic self-adjoint operator on an odd-dimensional manifold the
$\eta$-invariant of Atiyah--Patodi--Singer \cite{APS1}, which is by
definition equal to the value at the origin of the analytic continuation
of the $\eta$-function, is a finite number, i.e., the $\eta$-function does
not have a pole at the origin. Assuming additionally that $A$ is an
admissible operator of even order on an odd-dimensional manifold, one can
claim that the reduction of the $\eta$-invariant modulo ${\bf Z}$ is
invariant under deformations of the operator, while the integer jumps
occur as a result of (discontinuous) changes of the nonnegative spectral
subspace of the operator.

Such a homotopy invariance suggests that for this class of operators
the $\eta$-invariant of the operator $A$ is completely determined
by its nonnegative spectral subspace (this subspace is, actually,
the range of an even projection). This idea is realized in the following
proposition.

\begin{proposition}
\label{connect}Let $P_{+}\!\!$ be an even pseudodifferential
projection that is equal to a nonnegative spectral projection
for some admissible operator $A$. Then the dimension
of the projection with the trivial normalization
$\chi:K(M)\to {\bf R}$ and the
$\eta$-invariant of the operator $A$ are equal:
\[
d\left( P_{+}\right) \!=\!\frac{\eta \left(0,A\right) +\dim \ker A}2
\stackrel{{\rm def}}{=}{\eta }\left( A\right) .
\]
\end{proposition}

\noindent {\em Proof\/}. The orders of the elements of the group $K\left(
P^{*}M\right)/K(M)$ are powers of 2. Hence, (see Sec.~\ref{perv} and the
paper \cite{Gil7}) for sufficiently large $N$ the operator $2^NA$ is
homotopic (in the space of admissible self-adjoint elliptic operators) to
some operator, denoted by $A^{\prime },$ which is equal to a direct sum of
a positive and a negative admissible operator. Hence, we obtain ${\eta
}\left( A^{\prime }\right) =0$ (see \cite{Gil7}). Let us denote an
arbitrary homotopy of this form by $A_t.$ Recalling that the spectral flow
of the family $A_t$ through the point zero is equal to the net number of
jumps of the ${\eta }$-invariant, we obtain
\[
{\eta }\left( 2^NA\right) =-\limfunc{sf}A_t.
\]
On the other hand, by considering a spectral section
(see \cite{MePi1}) $P_t$ for the family $A_t$, we can show
(this is one of the definitions of the spectral flow) that
\[
\limfunc{sf}A_t=-\limfunc{ind}\left( 2^NP_{+},P_0\right) +\limfunc{ind}%
\left( P_{+,A^{\prime }},P_1\right) .
\]
The last expression is equal to the $d$-dimension (see Remark \ref{remmu})
\[
-\limfunc{ind}\left( 2^NP_{+},P_0\right) +\limfunc{ind}\left( P_{+,A^{\prime
}},P_1\right) =
-d(2^NP_+)+d(P_{+,A'})
=-d\left( 2^NP_{+}\right) .
\]
We finally obtain the desired formula
\[
{\eta }\left( A\right) =\frac 1{2^N}
{\eta }\left(
2^NA\right) =\frac 1{2^N}d\left( 2^NP_{+}\right) =d\left( P_{+}\right) .
\]

\begin{remark}
{\em An even projection $P$ satisfies the conditions of Proposition
\ref{connect} if and only if it is an admissible operator of order zero
itself.}
\end{remark}

Indeed, the admissibility of the operator
$A$ implies the admissibility of its spectral projection (see \cite{Gil7}).
To prove the converse statement, let us note that admissible operators
form an algebra. Consider an elliptic operator defined by the formula
\[
A=P\Delta P-\left( 1-P\right) \Delta \left( 1-P\right),
\]
where $\Delta $ denotes an arbitrary admissible positive self-adjoint
operator of order two with the principal symbol of the Laplacian $\left|
\xi \right| ^2.$ The operator $A$ is exactly the desired one (i.e. an
operator with nonnegative spectral projection equal to $P$). This
completes the proof of the remark.

In the paper \cite{Gil7}, the topological meaning of the invariant
${\eta }\,\func{mod}{\bf Z }$ (i.e. modulo jumps) was considered. It was
noticed that this reduction defines a homomorphism
\[
{\eta }:K\left( P^{*}M\right)/K(M) \longrightarrow {\bf Z}\left[
\frac 12\right] \func{mod}{\bf Z}.
\]
By virtue of Proposition~\ref{connect}, the
obtained index formula for operators in subspaces
(\ref{findx}), when reduced $\func{mod}{\bf Z,}$ gives a new
formula for the ${\eta }$-invariant,\footnote{But not for
all elements with principal symbols in the group
$K\left(
P^{*}M\right) .$}
which we state now as a corollary.

\begin{corollary}
\label{cor1}
Let
\begin{equation}
W\in {\rm ker}\left\{\pi^*:
K\left( P^{*}M\right)/K(M) \longrightarrow K\left( S^{*}M\right)/K(M)
\right\},\quad \pi :S^{*}M\rightarrow P^{*}M.  \label{requir1}
\end{equation}
Then the fractional part of the $\eta$-invariant of the bundle
$W$ is the half-integer
\begin{equation}
{\eta }\left( W\right) \equiv \frac12{\rm ind}_t
\left( \sigma(\xi)\sigma^{-1}(-\xi) :\pi_S^{*}F\rightarrow
\pi_S^{*}F\right) \func{mod}{\bf Z},\quad \pi_S:
S^{*}M\rightarrow M,F\in \limfunc{%
Vect}\left( M\right) ,  \label{usefull}
\end{equation}
where $\sigma $ is an arbitrary isomorphism of vector bundles
$\pi ^{*}W$ and $\pi_S^{*}F$ over the cotangent sphere bundle,
which exists by condition {\em (\ref{requir1}).}
\end{corollary}

The index formula for even boundary value problems leads to
a similar formula (cf. \cite{Gil8}).
\begin{corollary}
\label{cor2}
\emph{
\emph{Let }$\left( D,B,P\right) $\emph{\ be an even elliptic
boundary value problem on the manifold }$M$\emph{%
. Then the fractional part of the }$\eta$\emph{-invariant
of the even bundle }
$\func{Im}\sigma \left(
P\right) $\emph{\ is a half-integer and can be computed by the formula}
\[
\eta \left( \func{Im}\sigma \left( P\right) \right) \equiv \frac 12\limfunc{%
ind}_t\left( \sigma \left( D\right) \left( \xi \right) \cup \sigma%
\left( D\right) \left( -\xi \right) \right) \mbox{mod}{\bf{Z}}.
\]
}
\end{corollary}

The last statement gives the cobordism invariance, which we state now in
terms of $K$-theory. Let us recall some constructions from \cite{APS3}.

Consider the homomorphism
\begin{equation}
\begin{array}{ccc}
K\left( S^{*}X\right) /K\left( X\right)  & \longrightarrow  & K_c\left(
T^{*}X\oplus 1\right),  \\
W=\func{Im}P\subset \pi ^{*}E & \longrightarrow
 & \left[ \tau +i\left( 2P-1\right) \left|
\xi \right| \right],
\end{array}
 \label{bvpik}
\end{equation}
where $E\in \limfunc{Vect}\left( X\right) ,\pi :S^{*}X\rightarrow X
$ is the natural projection and the element
\[
\left[ \tau +i\left( 2P-1\right) \left| \xi \right| \right] \in K_c\left(
T^{*}X\oplus 1\right)
\]
is understood in the sense of difference construction, i.e. outside of the
zero section of the bundle $T^{*}X\oplus 1$ it defines an isomorphism in
the pullback of the bundle $E$ on the space $T^{*}X\oplus 1$. It can be
shown (see \cite{APS3}) that this homomorphism coincides with the
composition
\begin{equation}
K\left( S^{*}X\right) /K\left( X\right) \stackrel{\delta }{\rightarrow }%
K^1\left( T^{*}X\right) =K_c\left( T^{*}X\times {\mathbf R}\right) ,
\label{simi}
\end{equation}
where $\delta $ is induced by the coboundary operator in $K$-theory.
Moreover, the map $\delta $ is an isomorphism (for an odd-dimensional
manifold $X$ this is easy to show by using a nonvanishing vector field).
Making use of the identification (\ref{simi}) given by formula
(\ref{bvpik}), we obtain the desired cobordism invariance in the following
form.

\begin{corollary}
\label{cor3}
{\em
\emph{On the boundary }$%
X=\partial M$\emph{\ of a smooth even-dimensional manifold }$M$\emph{
\ consider an even subbundle }$\left[ W\right] \in
K\left( P^{*}X\right) /K\left( X\right) $ \emph{and
its pullback to the bundle of cotangent spheres,}
\[
\pi ^{*}\left[ W\right] \in K\left( S^{*}X\right) /K\left( X\right) ,\qquad
\pi :S^{*}X\rightarrow P^{*}X.
\]
\emph{Suppose that this extends inside the manifold $M$,
i.e. it lies in the range of the restriction operator}
\[
\begin{array}{c}
i^{*}:K_c\left( T^{*}M\right) \rightarrow K_c\left( \left. T^{*}M\right|
_X\right) ,\quad \left. T^{*}M\right| _X\subset T^{*}M, \\
\pi ^{*}\left[ W\right] \in \func{Im}i^{*}
\end{array}
\]
\emph{for the identification }(\ref{simi})\emph{.
Then the doubled }%
$\eta $\emph{-invariant of the bundle }$\;W$\emph{\
is an integer:}
\[
2\eta \left( W\right) {\rm mod}\, {\mathbf Z}=0.
\]
}
\end{corollary}

\section{Examples\label{xamples}}

\noindent
{\bf 1.} On a closed, connected, odd-dimensional, oriented, Riemannian manifold $M$
consider an elliptic self-adjoint differential operator of second order
acting on exterior $1$-forms by the formula
\[
A=d\delta -\delta d:\Lambda ^1\left( M\right) \longrightarrow \Lambda
^1\left( M\right) ,
\] where $d$ is the usual exterior derivative and $\delta $ is the adjoint
operator. A direct calculation shows that the principal symbol of the
nonnegative spectral projection $P$ for this operator at a point $\xi \neq 0$
is a projection on the line generated by the covector $\xi $ itself. In
other words, the projection $P$ is an even projection, while the image of
the principal symbol $\func{Im}\sigma ^{\prime }\left(P\right) $ on each
of the projective spaces $P_x^{*}M $ is exactly the tautological line
bundle, which is known to be nontrivial. At the same time, the pullback of
this bundle to the cosphere bundle is already trivial with the natural
trivialization
\begin{eqnarray*}
\kappa :\func{Im}\sigma \left( P\right) &\longrightarrow& {\bf C,} \\
\kappa \left( x,\xi ,\eta \right) &=&\left( \xi ,\eta \right) ,
\end{eqnarray*}
where $\left( \xi ,\eta \right) $ is the inner product
of two proportional vectors. The calculation made in the paper
\cite{Gil7} shows that in this case
\[
{\eta }\left( A\right)=\dim \ker \left. \left( d\delta +\delta
d\right) \right| _{\Lambda ^1\left( M\right) }-1=\dim H^1\left( M\right) -1.
\]
Hence, for an elliptic operator
\[
D:\func{Im}P\longrightarrow C^\infty \left( M\right)
\]
acting in the corresponding subspaces, the index formula can be rewritten
in the form
\[
\limfunc{ind}\left( D,P,P_{C^\infty \left( M\right) }\right) ={\rm ind}_t\left(
\sigma \left( D\right) \kappa ^{-1}:\pi^{*}{\bf C}\rightarrow \pi^{*}
{\bf C%
}\right) +\dim H^1\left( M\right) -1,
\]
where $$
{\rm ind}_t\left( \sigma \left( D\right) \kappa ^{-1}:\pi^{*}{\bf C}%
\rightarrow \pi^{*}{\bf C}\right)
$$
is the topological index of a usual elliptic pseudodifferential operator
acting on functions.\footnote{For ${\rm dim}M\ge 3$ this index is well known
to be zero. In order to obtain a nontrivial
index in this situation, one can consider matrix operators.}

\noindent
{\bf 2.} An operator similar to the one from the previous example is known
for the case of coefficients in a bundle $W\in \limfunc{Vect}\left(
M\right) $ and also for the case of forms of higher degrees. Namely,
consider a vector bundle $W$ with a connection $D.$ The corresponding
operator is defined by the formula
\[
A\otimes 1_W=DD^{*}-D^{*}D:\Lambda ^k\left( M,W\right) \longrightarrow
\Lambda ^k\left( M,W\right).
\]
It is an elliptic self-adjoint operator. Its principal symbol is
the tensor product
\[
\sigma \left( A\otimes 1_W\right) =\sigma \left( A\right) \otimes 1_W.
\]
In Gilkey's paper \cite{Gil7}, the problem of computation of the
fractional part of the $\eta$-invariant for operators of the above type
was posed. It follows from the formula (\ref{usefull}), for example, that
on a parallelizable manifold the fractional part of such operators acting
on 1-forms with coefficients in an arbitrary bundle $W$ is zero. Indeed,
by the index formula it is equal to half the index of an operator with a
constant principal symbol.

The problem posed by Gilkey asks for an operator with a nonzero fractional
part of the $\eta$-invariant. It remains open. In the context of even
projections this problem can be partly restated as follows: is it possible
that fractional terms appear in our index formulas?


\bibliographystyle{unsrt1}

\noindent
{\it Moscow State University}
\end{document}